\documentclass[letterpaper]{article}
\usepackage[T1]{fontenc}
\usepackage[latin1]{inputenc}
\usepackage{amsmath,amsthm,amsfonts,amssymb, euscript, stmaryrd, graphicx,superbox,appendix}
\usepackage{algorithm}
\usepackage[noend]{algpseudocode}
\usepackage{color}
\usepackage{pstricks}
\usepackage{psfrag}
\usepackage{amsthm}
\usepackage{fullpage}

\definecolor{colKeys}{rgb}{0,0,1}
\definecolor{unidentified}{rgb}{0,0,0}
\definecolor{concealments}{rgb}{0,0.5,0}
\definecolor{colString}{rgb}{0.6,0.1,0.1}

\definecolor{Grey1}{rgb}{0.5,0.5,0.5}

\def\myqed {{
\parfillskip=0pt 
\widowpenalty=10000 
\displaywidowpenalty=10000 
\finalhyphendemerits=0 
\leavevmode 
\unskip 
\nobreak 
\hfil 
\penalty50 
\hskip.2em 
\null 
\hfill 
$\square$
\par}} 

\newcommand{\footnoteremember}[2]{
 \footnote{#2}
 \newcounter{#1}
 \setcounter{#1}{\value{footnote}}
}
\newcommand{\footnoterecall}[1]{
 \footnotemark[\value{#1}]
}

\def\independent{{\perp\!\!\!\!\perp}}
\def\simdist{\stackrel{\mathcal{L}}{\sim}}

\author
{
Sylvain \textsc{Corlay}\footnote{Bloomberg L.P. Quantitative Research, 731 Lexington avenue, New York, NY 10022, USA. scorlay@bloomberg.net.}
\footnoteremember{myfootnote}{Laboratoire de Probabilités et Modèles Aléatoires, UMR 7599, Université Paris 6, case 188, 4, pl. Jussieu, F-75252 Paris Cedex 5, France. gilles.pages@upmc.fr}\\
 \and 
 Gilles \textsc{Pagès}\footnoterecall{myfootnote}
}

\title{Functional quantization-based stratified sampling methods}

\newtheorem{theo}{Theorem}[section]
\newtheorem{prop}[theo]{Proposition}

\newtheorem*{remark}{Remark}

\newtheorem{defi}{Definition}

\newcommand{\1}{\textbf{1}}
\newcommand{\N}{\mathbb{N}}

\newcommand{\R}{\mathbb{R}}

\newcommand{\E}{\mathbb{E}}

\newcommand{\PP}{\mathbb{P}}

\newcommand{\supp}{\operatorname{supp}}
\newcommand{\card}{\operatorname{card}}

\newcommand{\sspan}{\operatorname{span}}
\newcommand{\cov}{\operatorname{cov}}
\newcommand{\var}{\operatorname{Var}}
\newcommand{\Proj}{\operatorname{Proj}}

\newcommand{\ui}{{\underline{i}}}

\def\keywordname{{\bf Keywords:}} 
\newcommand{\keywords}[1]{\par\addvspace\baselineskip\noindent\keywordname\enspace\ignorespaces#1}

\graphicspath{{./}{figures/}}

\begin{document}

\maketitle

\begin{abstract}
\par In this article, we propose several quantization-based stratified sampling methods to reduce the variance of a Monte Carlo simulation.
\par Theoretical aspects of stratification lead to a strong link between optimal quadratic quantization and the variance reduction that can be achieved with stratified sampling. We first put the emphasis on the consistency of quantization for partitioning the state space in stratified sampling methods in both finite and infinite dimensional cases. We show that the proposed quantization-based strata design has uniform efficiency among the class of Lipschitz continuous functionals. 
\par Then a stratified sampling algorithm based on product functional quantization is proposed for path-dependent functionals of multi-factor diffusions. The method is also available for other Gaussian processes such as Brownian bridge or Ornstein-Uhlenbeck processes. We derive in detail the case of Ornstein-Uhlenbeck processes. 
\par We also study the balance between the algorithmic complexity of the simulation and the variance reduction factor. 
\end{abstract}
\keywords{functional quantization, vector quantization, stratification, variance reduction, Monte Carlo simulation, Karhunen-Loève, Gaussian process, Brownian motion, Brownian bridge, Ornstein-Uhlenbeck process, Ornstein-Uhlenbeck bridge, principal component analysis, numerical integration, option pricing, Voronoi diagram, product quantizer, path-dependent option.} 

\section*{Introduction}\label{sec:introduction}
\par The quantization of a random variable $X$ consists of its approximation by a random variable $Y$ taking finitely many values. This problem has been initially investigated for its applications to signal transmission and for compression issues \cite{GershoGrayVectorQuantization}. In this context, quantization is a method of signal discretization. The aim is to choose the random variable $Y$ so as to minimize the resulting error for a fixed quantization level $N$.
\par More recently, quantization was introduced in numerical probability to devise numerical integration methods \cite{PagesGaussianQuantization} and to solve multidimensional stochastic control problems such as the pricing of American options \cite{BallyPagesPrintemsAmerican1} and swing options \cite{BardouBouthemyPagesSwing1}. 
Optimal quantization has many other applications and extensions in various fields such as automatic clustering (quantization of empirical measures) and pattern recognition. 
\par Since the early $2000$'s, the infinite-dimensional setting has been extensively investigated from both theoretical and numerical viewpoints with a special attention paid to functional quantization \cite{LuschgyPagesFunctional3,PagesPrintemsFunctional4}. Bi-measurable stochastic processes are viewed as random variables valued in functional spaces. 
\par Still Monte Carlo simulations remain the most common numerical method in the field of numerical probability. One reason is that it is easy to implement in an industrial configuration. In the industry of derivatives, banks implement generic Monte Carlo frameworks for pricing and hedging their positions with a wide variety of financial products and models. Besides, Monte Carlo simulations are easily parallelized. 
\par Variance reduction methods can be used to dramatically reduce the computation time of a Monte Carlo simulation, or to increase its accuracy. The main variance reduction methods are (adaptive) control variate, pre-conditioning, importance sampling and stratification \cite{GlassermanMonteCarlo,LemairePagesSampling}. The problem is that these methods may strongly depend on the payoff or the model and require significant changes in the practical implementation of the Monte Carlo simulation. Therefore, most practitioners do not use the most sophisticated methods except for marginal cases. 
\par In this article, we point out theoretical aspects of quantization that idraw a strong link between the problem of optimal quadratic quantization of a random variable and the variance reduction that can be achieved by stratification. We emphasize the consistency of quantization for designing strata in stratified sampling methods in both the finite and infinite dimensional settings. Then we devise a stratified sampling algorithm based on product functional quantization for path-dependent functionals of multi-factor Brownian diffusions. We show that this strata design has uniform efficiency among the class of Lipschitz continuous functionals of Brownian motion. The simulation cost of the conditional path is $O(n)$ where $n$ is the number of discretization dates as in the naive unconditioned Monte Carlo simulations. In this context, the proposed approach can be considered as a guided Monte Carlo simulation (see Figure \ref{fig:brownian_motion_conditional_simulation}). 
The method is applicable with any Gaussian process as soon as we can derive its Karhunen-Loève expansion. This is the case for Brownian bridge and Ornstein-Uhlenbeck processes. The special case of Ornstein-Uhlenbeck processes is detailed in Appendix \ref{sec:karhunen_loeve_basis}. The case of the Ornstein-Uhlenbeck bridge is presented in \cite{CorlayOUBridge}.
\par A very common situation is the case of Monte Carlo simulations of multi-factor Brownian diffusions approximated with an Euler discretization or another time-discretization scheme. The presented method is particularly well suited for this case, regardless of how the Brownian paths are used in the model, to drive the dynamics of the stock price, a volatility process or a drift term. Functional stratification can be used as a generic variance reduction method which does not require a reimplementation of the whole framework but only the way it is input with Brownian paths.
\par The article is organized as follows. Section \ref{sec:quantization_base} presents some necessary background on optimal quantization. The emphasis is on the functional quantization of Gaussian processes. Section \ref{seq:quantization_variance_reduction} briefly covers the first functional quantization-based variance reduction method that was proposed in \cite{PagesPrintemsFunctional4,LejayReutenauerControl}. Section \ref{sec:stratification} outlines the links between quantization and stratification with an emphasis on the Gaussian case. The method is further detailed in the infinite-dimensional case for Gaussian processes in Section \ref{sec:functional_stratification}. We present a simulation method for the case of Brownian motion and other examples of Gaussian processes (such as Brownian bridge and Ornstein-Uhlenbeck processes) that preserves the $O(n)$ simulation complexity where $n$ is the number of time steps. In Section \ref{sec:option_pricing}, we provide numerical experiments of the method with option pricing problems arising in mathematical finance. Appendix \ref{sec:karhunen_loeve_basis} presents the computation of the Karhunen-Loève expansion of Ornstein-Uhlenbeck processes. Appendix \ref{sec:annex_R2_computation} presents the derivation of closed-form expressions of some regression matrices needed for our stratified sampling algorithm, in the cases of Brownian motion, Brownian bridge and Ornstein-Uhlenbeck processes.
\section{Vector and functional quantization}\label{sec:quantization_base}
\subsection{Introduction to quantization of random variables}
\par Let $(\Omega,\mathcal{A},\PP)$ be a probability space and $(E, |\cdot|)$ a reflexive separable Banach space. The principle of the quantization of a random variable $X$ taking its values in $E$ is to approximate $X$ by a random variable $Y$ taking a finite number $N$ of values in $E$. The discrete random variable $Y$ is a quantizer of $X$ of level $N$. The resulting discretization error to be minimized is the $L^p$ norm of $|X-Y|$. 
\begin{equation}\label{eq:minimization_problem}
\min \left\{ \| X-Y\|_p, Y:\Omega \to E \textnormal{ measurable}, \ \card(Y(\Omega)) \leq N \right\}. 
\end{equation}
\begin{defi}[Voronoi partition] Consider $N \in \N^*$, $\Gamma = \{\gamma_1,\ldots,\gamma_N\} \subset E$ and let $C = \{C_1,\ldots,C_N\}$ be a Borel partition of $E$. $C$ is a Voronoi partition associated with $\Gamma$ if $\forall i \in \{ 1,\ldots, N \}, \ C_i \subset \{\xi \in E, |\xi -\gamma_i | = \min\limits_{j \in \{ 1, \ldots, N \} } | \xi - \gamma_j | \}$. $C_i$ is called Voronoi cell associated with $\gamma_i$ in $C$.
\end{defi}. 
\begin{prop}\label{thm:best_nnp}
	\par Let $X$ and $Y$ be two random variables valued in $E$, where $Y$ takes its values in the fixed set of knots $\Gamma = \{\gamma_1, \ldots, \gamma_N\} \subset E$ for $N\in \N^*$. We define $\widehat{X}^\Gamma := \Proj_{\Gamma}(X)$ where $\Proj_{\Gamma} = \sum\limits_{i=1}^N \gamma_i \1_{C_i}$ is a nearest neighbor projection onto $\Gamma$. Then we have $\left|X-\widehat{X}^\Gamma\right| \leq |X-Y| \textnormal{ a.s.}$ and thus $\left\| X-\widehat{X}^\Gamma \right\|_p \leq \| X-Y\|_p$.
\end{prop}
\par \noindent A consequence is that solving \eqref{eq:minimization_problem} amounts to solving the simpler problem
$$
\min \left\{ \| X-\Proj_{\Gamma}(X) \|_p, \ \Gamma \subset E, \ \card(\Gamma) \leq N\right\}.
$$
\par \noindent The quantity $\| X-\Proj_{\Gamma}(X) \|_p$ is called the mean $L^p$ quantization error. The problem of the existence of a minimum is addressed in \cite{PagesIntegVectorQuant,GrafLushgyMonograf} for the finite-dimensional case. 
	\begin{itemize}
		\item For every $N \geq 1$, the mean $L^p$ quantization error is Lipschitz continuous and reaches a minimum. An $N$-tuple that achieves the minimum has pairwise distinct components, as soon as $\card(\supp(\PP_X)) \geq N$. 
    \end{itemize}		
\par \noindent This result stands in the general case of a random variable valued in a reflexive separable Banach space \cite{LuschgyPagesFunctional3}.
    \begin{itemize}
		\item If the support of $\PP_X$ has infinite cardinal, the optimal quantization error decreases, and converges to $0$ as the quantization level $N$ goes to infinity. In the finite-dimensional case, and for distributions that are absolutely continuous with respect to the Lebesgue measure, the rate of convergence is ruled by Theorem \ref{thm:buckley_wise_graf_luschgy}. 
	\end{itemize}

\begin{theo}[Zador, Luschgy, Pagès]\label{thm:buckley_wise_graf_luschgy}
	\begin{itemize}
	\item \emph{(Sharp rate)} Let $\ r > 0$ and $X:\Omega \to \R^d \in L^{r+\eta}(\PP)$ for some $\eta>0$. Let $\PP_X(d\xi) = \phi(\xi) d \xi + \mu(d\xi)$ be the canonical decomposition of $\PP_X$ ($\mu$ and the Lebesgue measure are singular). Then, if $\phi \not\equiv 0$, the $L^r$ quantization error at level $N$, $\mathcal{E}_{N,r}$ satisfies
	\begin{equation}
	\mathcal{E}_{N,r}(X,\R^d) \underset{N \to \infty}{\sim} \widetilde{J}_{r,d} \times \left(\int_{\R^d} \phi^{\frac{d}{d+r}}(u) du \right)^{\frac{1}{d}+\frac{1}{r}} \times N^{-\frac{1}{d}}, \quad \textnormal{where} \quad \widetilde{J}_{r,d} \in (0,\infty).
        \end{equation}
	\item \emph{(Non-asymptotic upper bound)} Let $d \geq 1$. There exists $C_{d,r,\eta} \in (0,\infty)$ such that, for every $\R^d$-valued random vector $X$,
	\begin{equation}\label{eq:non_asymp}
	\forall N \geq 1, \hspace{5mm} \mathcal{E}_{N,r}(X,\R^d) \leq C_{d,r,\eta}\|X\|_{r+\eta} N^{-\frac{1}{d}}.
	\end{equation}
	\end{itemize}
\end{theo}
\par \noindent The first claim was stated for the case of distributions with compact support by Zador in \cite{ZadorAsymptoticError}. The extension to general probability distributions in $\R^d$ was developed in \cite{BuckleyWise}. The first mathematically rigorous proof can be found in \cite{GrafLushgyMonograf}. The non-asymptotic error bound \eqref{eq:non_asymp} is proved in \cite{FunctionalQuantizationLevy}. 
\par
\begin{figure}[!ht]
 	\begin{minipage}[c]{.46\linewidth}
 	\psfrag{-4}{$-4$}
	\psfrag{-3}{$-3$}
	\psfrag{-2}{$-2$}
	\psfrag{-1}{$-1$}
	\psfrag{0}{$0$}
	\psfrag{1}{$1$}
	\psfrag{2}{$2$}
	\psfrag{3}{$3$}
	\psfrag{4}{$4$}
	\includegraphics[height=5.5cm]{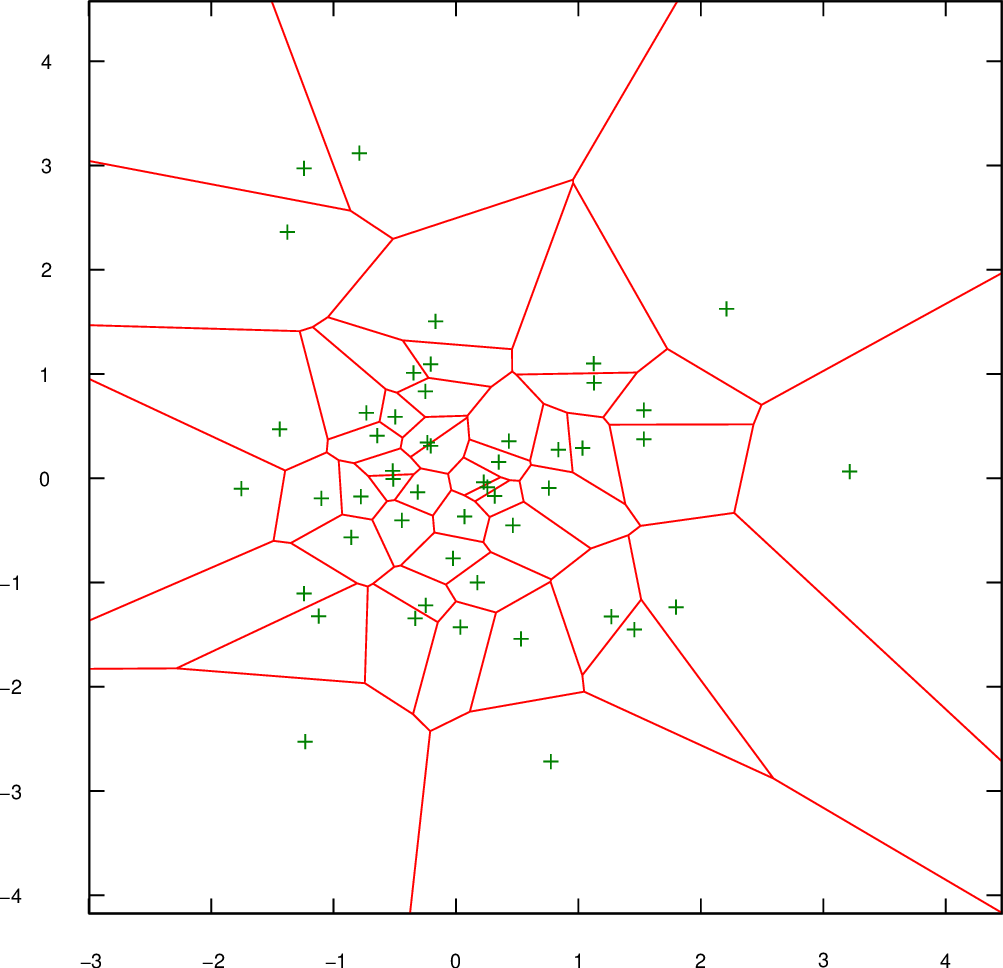}
	\end{minipage} \hfill
	\begin{minipage}[c]{.46\linewidth}
	\psfrag{-3}{$-3$}
	\psfrag{-2}{$-2$}
	\psfrag{-1}{$-1$}
	\psfrag{0}{$0$}
	\psfrag{1}{$1$}
	\psfrag{2}{$2$}
	\psfrag{3}{$3$}
	\includegraphics[height=5.5cm]{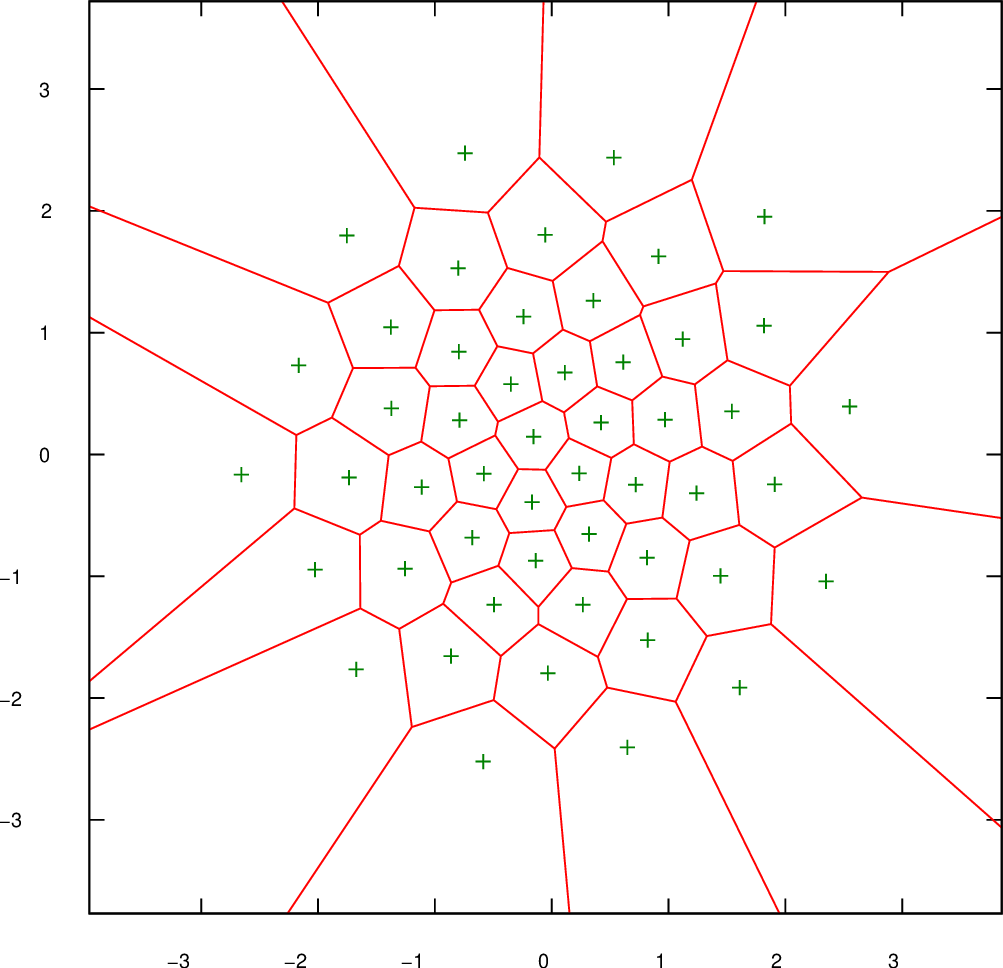}
	\end{minipage}
 	\caption{Voronoi partition of a random quantizer (left) and an optimal quantizer (right) of level $N=48$ of the $\mathcal{N}(0,I_2)$ distribution.}
	\label{fig:gaussian_voronoi}
\end{figure}

 In Figure \ref{fig:gaussian_voronoi}, we display the Voronoi partition of a random $N$-quantizer and an optimal quadratic quantizer of level $N$ for the bivariate normal distribution $\mathcal{N}(0,I_2)$. 
\subsection{Self-consistency of optimal quantizers}\label{sec:stationarity}
\par We now assume that $E$ is a separable Hilbert space $(H,\langle\cdot,\cdot\rangle_H)$. We denote by $\mathcal{C}_N(X)$ the set of $L^2$ optimal quantizers of $X$ of level $N$, and by $\mathcal{E}_N(X)$ the minimal quadratic distortion that can be achieved when approximating $X$ by a quantizer of level $N$.
\begin{defi}[Stationarity]
\par A quantizer $Y$ of $X$ is stationary (or self-consistent) if
\begin{equation}\label{eq:stationarity}
Y = \E[X | Y].
\end{equation}
\end{defi}
\begin{prop}[Stationarity of $L^2$ optimal quantizers]\label{prop:stationarity}
\par A quadratic optimal quantizer is stationary. 
\end{prop}
\par \noindent We refer to \cite{GrafLushgyMonograf} for the proof in the finite-dimensional setting and to \cite{LuschgyPagesFunctional3} for the more general case of separable Hilbert spaces. Stationarity is a particularity of the quadratic case ($p=2$). In the general $L^p$ setting, a similar property involving the notion of $p$-center holds \cite{GLPApprox}.
\begin{prop}
\par Let $X$ be an $H$-valued $L^2$ random variable. Let us denote by $D_N^X$ the squared quadratic quantization error associated with a codebook of size $N$ with respect to $X$. 
$$
\begin{array}{lccc}
D_N^X :	& H^N &\to & \R_+ \\
		& \Gamma = (\gamma_1, \ldots, \gamma_N) & \mapsto & \E\big[ \min\limits_{1 \leq i \leq N} |X-\gamma_i|_H^2 \big].
\end{array}
$$
\par \noindent The distortion $D_N^X$ is $|\cdot|_H$-differentiable at $N$-quantizers $\Gamma \in H^N$ with pairwise distinct components and such that boundaries of Voronoi cells are $\PP_X$-negligible
\begin{equation}\label{eq:stationarity_critical}
\nabla D_N^X(\Gamma) = 2 \Big( \int_{C_i(\Gamma)} (\gamma_i - \xi) \PP_X(d \xi) \Big)_{1 \leq i \leq N} = 2 \left( \E \left[ \left(\widehat{X}^{\Gamma} - X \right) \1_{\left\{\widehat{X}^{\Gamma} = \gamma_i\right\}} \right] \right)_{1 \leq i \leq N}.
\end{equation}
Hence any Voronoi quantizer associated with a critical point of $D_N^X$ is a stationary quantizer.
\end{prop}
\par \noindent We refer to \cite{PagesReview} for a detailed proof. 
\begin{defi}[Centroidal projection]\label{def:centroidal_projection}
	\par Let $C = \{C_1, \ldots, C_N\}$ be a Borel partition of $H$. For $1\leq i\leq N$, we define
 	$G_i := 	\left\{\begin{array}{llll}
 			\E[X | X \in C_i] & \textnormal{if } \PP[X \in C_i] \neq 0,\\
 			0 & \textnormal{otherwise,}
 			\end{array}\right.$
 	the centroids associated with $X$ and $C$. 
 	\par \noindent The centroidal projection associated with $C$ and $X$ is the application $\Proj_{C,X} : x \mapsto \sum\limits_{i=1}^N G_i \1_{C_i}(x)$.
\end{defi}

\subsection{Optimal quantization and principal component analysis}\label{sec:quantization_pca}
For any finite-dimensional subspace $U$ of $H$, we denote by $\Pi_U$ the orthogonal projection onto $U$. 
\begin{prop}\label{prop:linear_subspace}
\par Let $U$ be a finite-dimensional linear subspace of $H$. Then 
\begin{equation}\label{eq:decomposition_quant_pca}
\mathcal{E}_N(X)^2 \leq \E \left[ \left|X-\Pi_U(X)\right|^2\right] + \mathcal{E}_N(\Pi_U(X))^2.
\end{equation}
Moreover, if an optimal quantizer of $X$ of level $N$ lies in $U$, we have equality in \eqref{eq:decomposition_quant_pca}.
\end{prop}
We refer to \cite{LuschgyPagesFunctional3} for a detailed proof. This allows us to define the quantization dimension of $X$ of level $N$ by $d_N(X) := \min \big\{ \dim \sspan(\Gamma), \Gamma \in \mathcal{C}_N(X) \big\}$. It follows from Proposition \ref{prop:linear_subspace} that 
$$
\mathcal{E}_N^2(X) = \min \left\{ \E[\|X-\Pi_V(X) \|^2] + \mathcal{E}_N^2(\Pi_V(X)), \begin{array}{lll}V \subset H \textnormal{ linear subspace }\\ \textnormal{such that} \dim V \geq d_N(X)\end{array} \right\}.
$$
\subsubsection{Covariance operator of a Gaussian measure}
\begin{defi}
\par Let $X$ be a centered $H$-valued $L^2$ Gaussian random variable. Its covariance operator $C_X: H \to H$ is defined by $C_X y = \E[ \langle y, X \rangle X ]$.
\begin{enumerate}
\item If $X$ is $\R^d$-valued, the matrix of $C_X$ in the canonical basis is the covariance matrix of $X$.
\item If $X=(X_t)_{t \in [0,T]}$ is a bi-measurable centered process of covariance function $\Gamma_X(s,t) := \E[X_s X_t]$ satisfying $\int_{[0,T]} \Gamma_X(s,s) ds < +\infty$, then $X$ can be seen as a random variable valued in $L^2([0,T],dt)$ satisfying $\E\left[\left|X\right|^2\right] < \infty$, and 
$$
C_X y = \int_{[0,T]} y(s) \Gamma_X(s,\cdot)ds, \quad y \in L^2([0,T],dt).
$$
\end{enumerate}
\end{defi}
\par \noindent In \cite{LuschgyPagesFunctional3}, it is proved that linear subspaces $U$ of $H$ spanned by $n$-stationary quantizers of Gaussian measures correspond to principal subspaces of $X$. In other words, they are spanned by the eigenvectors of $C_X$ corresponding to the largest eigenvalues.  
\begin{theo}
Let $\Gamma$ be an optimal codebook for the Gaussian random variable $X$, $U = \sspan(\Gamma)$ and $m= \dim U$. Then $C_X(U) = U$ and
$\E \left[\left|X-\Pi_U(X)\right|^2\right] = \sum\limits_{j \geq m+1} \lambda^X_j$, where $\lambda^X_1 \geq \lambda^X_2 \geq \cdots > 0$ are the ordered non-zero eigenvalues of $C_X$ (repeated as many times as their multiplicity). We have
\end{theo}
$$
\sum\limits_{j \geq m+1} \lambda^X_j = \inf \left\{ \E\left[\left|X-\Pi_V(X)\right|^2\right], V \subset H \textnormal{ linear subspace}, \ \dim V= m \right\}.
$$
\par \noindent The minimal quadratic distortion $\mathcal{E}_N(X)$ is given by
\begin{equation}\label{eq:gauss_infinite_finite1}
\mathcal{E}_N(X)^2 = \sum\limits_{j \geq m+1} \lambda^X_j + \mathcal{E}_N \left( \bigotimes\limits_{j=1}^m \mathcal{N}\left(0,\lambda^X_j\right) \right)^2 \ \textnormal{ for } m \geq d_N(X),
\end{equation}
\par \noindent A proof is available in \cite{LuschgyPagesFunctional3}. This shows that the optimal quantization of a Gaussian process $X$ boils down to a finite-dimensional quantization problem, if the Karhunen-Loève eigensystem $(e^X_n, \lambda^X_n)_{n \in \N^*}$ is known. 
\subsection{Product quantization}
\par Let $(e_n)_{n \in \N^*}$ be a Hilbert basis of $H$, and $(N_n)_{n \geq 1}$ an integer sequence such that $\prod_{n \geq 1} N_n < \infty$ (so that $N_n=1$ for large enough $n$). For every $n \in \N^*$, we consider a codebook of size $N_n$, $\Gamma^n := \left\{ \gamma_1^n, \ldots, \gamma_{N_n}^n \right\} \subset \R$. 
\par The codebook $\Gamma$ is defined as the set of knots in $H$ whose coordinates in the base $(e_n)_{n \in \N^*}$ are the Cartesian product of the one-dimensional codebooks $\Gamma^n$. 
\begin{prop}[Case of independent marginals]
\par With the same notation, if we assume that the marginals of $X$, $(\langle X,e_1 \rangle,\langle X, e_2 \rangle, \ldots)$ are independent, and that for each $k \in \N^*$,  $Y^k := \Proj_{\Gamma^k}(\langle X, e_k \rangle)$ is a stationary quantizer of $\langle X, e_k \rangle$, then $Y = \Proj_\Gamma(X)$ is a stationary quantizer of $X$. 
\end{prop}
\par \noindent In the case of independent marginals, optimal product quantization remains stationary and the simple shape of Voronoi cells simplifies the nearest neighbor search. 
\subsection{Numerical optimal quantization}
\par Various algorithms have been developed to compute optimal $N$-grids in the finite-dimensional setting. A review of these methods is available in \cite{PagesReview}. Let us mention Lloyd's algorithm for the quadratic case. Another approach is the stochastic gradient method which is suggested by the fact that the quadratic distortion function has an integral representation and is differentiable at any $N$-tuple having pairwise distinct components and a $\PP_X$-negligible Voronoi tessellation boundary \cite{PagesGaussianQuantization}.
\par Equation \eqref{eq:stationarity_critical} shows that any Voronoi quantizer associated with a critical point of $D_N^X$ is a stationary quantizer. In the case of one-dimensional distributions, such as the Gaussian distribution, the (tridiagonal) Hessian of the distortion has a closed-form expression. Hence, a Newton-Raphson method can be easily implemented. It is thoroughly studied in \cite{PagesGaussianQuantization} in the Gaussian case and remains the fastest way to compute $L^2$ optimal quantizers of one-dimensional Gaussian variables.
\subsection{Quantization of Gaussian processes}
\subsubsection{Optimal quantization}
\par From now on, we will assume that $X$ is a bi-measurable Gaussian process and has a continuous covariance function $\Gamma^X$ and satisfies $\E\left[|X|^2_{L^2_T}\right] = \int\limits_0^T \E[X^2_s]ds < \infty$. 
\par We have seen in Section \ref{sec:quantization_pca} that in this setting, the $L^2$ optimal quantization $X$ amounts to the quantization of a finite-dimensional Gaussian vector $\bigotimes\limits_{j=1}^m \mathcal{N}\left(0,\lambda^X_j\right)$ for some positive integer $m$, the quantization dimension.\par Several usual Gaussian processes have explicit Karhunen-Loève expansions, such as Brownian motion, Brownian bridge and Ornstein-Uhlenbeck processes and bridges. (The case of a stationary Ornstein-Uhlenbeck process is derived for normalized parameters in the stationary case in \cite[p.195]{HirschLacombe}.) In Appendix \ref{sec:karhunen_loeve_basis}, we derive the Karhunen-Loève expansion of the Ornstein-Uhlenbeck process in the general case (for any value of the parameters and the initial variance). The K-L expansion of the Ornstein-Uhlenbeck bridge is derived in \cite{CorlayOUBridge}. To the best of our knowedge, no closed-form expression is available for fractional Brownian motion. In the article, numerical examples will be presented for the following cases.
{
\small
\begin{enumerate}
	\item \textbf{Brownian motion on $[0,T]$}:
	\begin{equation}\label{eq:brownian_motion_kl}
	e_n^W(t) := \sqrt{\frac{2}{T}} \sin \left( \pi (n-1/2) \frac{t}{T}\right), \hspace{8mm} \lambda_n^W:= \left(\frac{T}{\pi (n-1/2)} \right)^2, \hspace{4mm} n \geq 1.
	\end{equation}
	\item \textbf{Brownian bridge on $[0,T]$}: 
	\begin{equation}\label{eq:brownian_bridge_kl}
	e_n^B(t) := \sqrt{\frac{2}{T}} \sin \left( \pi n \frac{t}{T}\right), \hspace{8mm} \lambda_n^B:= \left(\frac{T}{\pi n} \right)^2, \hspace{4mm} n \geq 1.
	\end{equation}
	\item \textbf{The Ornstein-Uhlenbeck process on $[0,T]$, starting from $0$}, and defined by the SDE
	\begin{equation}\label{eq:ornstein_uhlenbeck_diffusion}
	dX_t = -\theta X_t dt + \sigma dW_t,
	\end{equation}
	with $\sigma \geq 0$, $\theta>0$ and $W$ a standard Brownian motion on $[0,T]$:
	\begin{equation}
	e_n^{OU}(t) := \frac{1}{\sqrt{\frac{T}{2}-\frac{\sin(2\omega_n T )}{4\omega_n}}} \sin(\omega_n t), \hspace{8mm} \lambda_n^{OU}:= \frac{\sigma^2}{\omega_n^2 + \theta^2}, \hspace{4mm} n \geq 1,
	\end{equation}
	where $(\omega_n)_{n \geq 1}$ are the increasingly sorted positive solutions of $\theta \sin(\omega_n T) + \omega_n \cos(\omega_n T) = 0$ (see Appendix \ref{sec:karhunen_loeve_basis}).
	\item \textbf{The stationary Ornstein-Uhlenbeck process on $[0,T]$} (see Appendix \ref{sec:karhunen_loeve_basis}).
\end{enumerate}
}

\par
\begin{figure}[!ht]
	\begin{center}
	\psfrag{0}{$0$}
	\psfrag{0.2}{$0.2$}
	\psfrag{0.4}{$0.4$}
	\psfrag{0.6}{$0.6$}
	\psfrag{0.8}{$0.8$}
	\psfrag{1}{$1$}
	\psfrag{-2.5}{$-2.5$}
	\psfrag{-2}{$-2$}
	\psfrag{-1.5}{$-1.5$}
	\psfrag{-1}{$-1$}
	\psfrag{-0.5}{$-0.5$}
	\psfrag{0.5}{$0.5$}
	\psfrag{1.5}{$1.5$}
	\psfrag{2}{$2$}
	\psfrag{2.5}{$2.5$}
	\includegraphics[width=6.5cm]{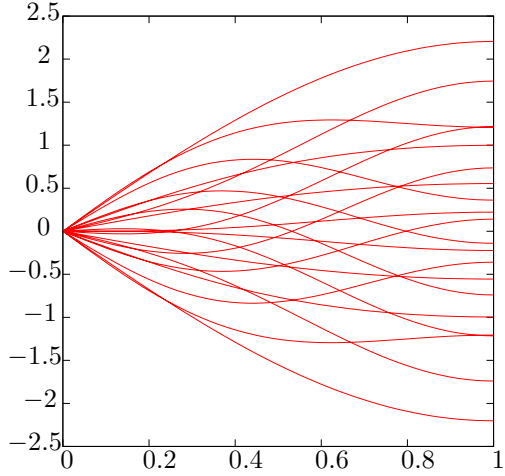}
	\caption{Optimal quantization of Brownian motion on $[0,1]$. }
	\label{fig:brownian_motion_optimal_quantizer}
	\end{center}
\end{figure}

 \noindent In Figure \ref{fig:brownian_motion_optimal_quantizer}, we display a quadratic $N$-optimal quantizer of Brownian motion. 
\subsubsection{Product quantization}\label{sec:functional_product_quantization}
\par Thanks to Equation \eqref{eq:gauss_infinite_finite1}, the product quantization of the finite-dimensional distribution $ \xi \simdist \bigotimes\limits_{j=1}^m \mathcal{N}\left(0,\lambda^X_j\right)$ yields a stationary quantizer $\widehat{X}$ of $X$ of the form $\widehat{X} = \sum\limits_{n \geq 1} \sqrt{\lambda^X_n} \widehat{\xi}_n e_n^X$, where $\widehat{\xi}_n$ is an optimal $N_n$-quantizer of $\xi_n$ and $\prod\limits_{n \geq 1} N_n \leq N, \ N_n \geq 1$ (so that for large enough $n$, $N_n = 1$ and $\widehat{\xi}_n \equiv 0$.) The paths corresponding to a multi-index $\ui = \{i_1, \ldots, i_n, \ldots \}$ are of the form 
$\chi_{\ui} = \sum\limits_{n \geq 1} \sqrt{\lambda_n^X} \gamma_{i_n}^{(N_n)} e_n^X.$
\par Such a functional quantizer $\widehat{X}$ is called a K-L product quantizer. Furthermore, we denote by 
$\mathcal{O}_{pq}(X,N)$ the set of K-L product quantizers of size at most $N$ of $X$. In the case of product quantization, the counterpart of Equation \eqref{eq:gauss_infinite_finite1} is 
\begin{multline}\label{eq:product_distortion}
\E\left[\min\limits_{\ui} \left|X-\chi_{\ui}\right|^2\right] = \sum\limits_{n=1}^m \lambda^X_n \E\left[ \min\limits_{1 \leq i_n \leq N_n} \left|\xi_n - \gamma_{i_n}^{(N_n)}\right|^2 \right] + \sum\limits_{n \geq m+1} \lambda^X_n \\
					= \sum\limits_{n=1}^m \lambda^X_n \left(\E\left[ \min\limits_{1 \leq i_n \leq N_n} \left|\xi_n - \gamma_{i_n}^{(N_n)}\right|^2 \right] -1 \right)+ \E\left[|X|^2_{L^2_T}\right],
\end{multline}
where $m$ is the quantization dimension. 
\subsubsection{Product decomposition and blind optimization}\label{sec:blind_optimization_quadratic}
\par The minimal quadratic error for a K-L product quantizer of level $N$ is the solution of the minimization problem
\begin{equation}\label{eq:blind_optimization_dist}
\mathcal{E}^{pq}_N := \min \Big\{\mathcal{E}(\chi), \ \chi \in \mathcal{O}_{pq}(X,N) \Big\},
\end{equation}
where $\mathcal{E}(\chi)$ is the quadratic distortion of the product quantizer $\chi$. Thanks to \eqref{eq:product_distortion}, this comes to
\begin{equation}\label{eq:blind_optimization_dist2}
\min \Big\{ \sum\limits_{n=1}^d \lambda^X_n \mathcal{E}_{N_n}\left(\mathcal{N}(0,1)\right)^2 + \sum\limits_{n \geq d+1} \lambda^X_n, \ N_1 \times \cdots \times N_d \leq N, \ d \geq 1 \Big\}. 
\end{equation}
A solution of \eqref{eq:blind_optimization_dist} is called an optimal K-L product quantizer. 
\par The blind optimization procedure consists of computing the criterion for every possible decomposition $N_1 \times \cdots \times N_d \leq N$, $d \geq 1$ and $N_1 \geq N_2 \geq \cdots$. For a given Gaussian process $X$, results can be kept off-line for a future use. The method is more thoroughly described in \cite{PagesPrintemsFunctional4}. Optimal decompositions for a wide range of values of $N$ for both Brownian bridge and Brownian motion are available on the web site {\verb www.quantize.maths-fi.com } \cite{WebSiteGaussian} for download. In the case of Ornstein-Uhlenbeck processes, the optimal decomposition depends on the diffusion parameters ($\sigma$ and $\theta$ in \eqref{eq:ornstein_uhlenbeck_diffusion}) and the maturity. 
\par Some optimal decompositions for the stationary Ornstein-Uhlenbeck process are given in Table \ref{tab:ornstein_uhlenbeck_optimal_decomposition}. 
\begin{table}[!ht]
\begin{center}
\begin{tabular}{|c|c|c|c|}
\hline
$N$ & $N_{rec}$ & Squared $L^2$ quantization Error & Product decomposition\\ 
\hline \hline
$1$ & $1$ & $1.5$ & $1$\\
$10$ & $10$	& $0.65318$ & $5 \ \times  \ 2$\\
$100$ & $96$ & $0.40929$ & $6 \ \times \ 4 \ \times \ 2 \ \times \ 2$\\
$1000$ & $960$ & $0.29618$ & $10 \ \times \ 6 \ \times \ 4 \ \times \ 2 \ \times \ 2$\\
$10000$ & $9984$ & $0.23150$ & $13 \ \times \ 8 \ \times \ 4 \ \times \ 3 \ \times \ 2 \ \times \ 2 \ \times \ 2$\\
\hline
\end{tabular}
\caption{Record of optimal product decompositions of the stationary centered Ornstein-Uhlenbeck process solution of the SDE $dX_t = - X_t dt+ dW_t$ on $[0,3]$. }
\label{tab:ornstein_uhlenbeck_optimal_decomposition}
\end{center}
\end{table}
\par In the following, we will face similar cases (other criteria than the quadratic distortion) where the blind optimization procedure applies. 
\vspace{1mm}
\par
\begin{figure}[!ht]
\begin{minipage}[c]{.46\linewidth}
\psfrag{0}{$0$}
\psfrag{0.2}{$0.2$}
\psfrag{0.4}{$0.4$}
\psfrag{0.6}{$0.6$}
\psfrag{0.8}{$0.8$}
\psfrag{1}{$1$}
\psfrag{-2.5}{$-2.5$}
\psfrag{-2}{$-2$}
\psfrag{-1.5}{$-1.5$}
\psfrag{-1}{$-1$}
\psfrag{-0.5}{$-0.5$}
\psfrag{0.5}{$0.5$}
\psfrag{1.5}{$1.5$}
\psfrag{2}{$2$}
\psfrag{2.5}{$2.5$}
\includegraphics[height=5.5cm]{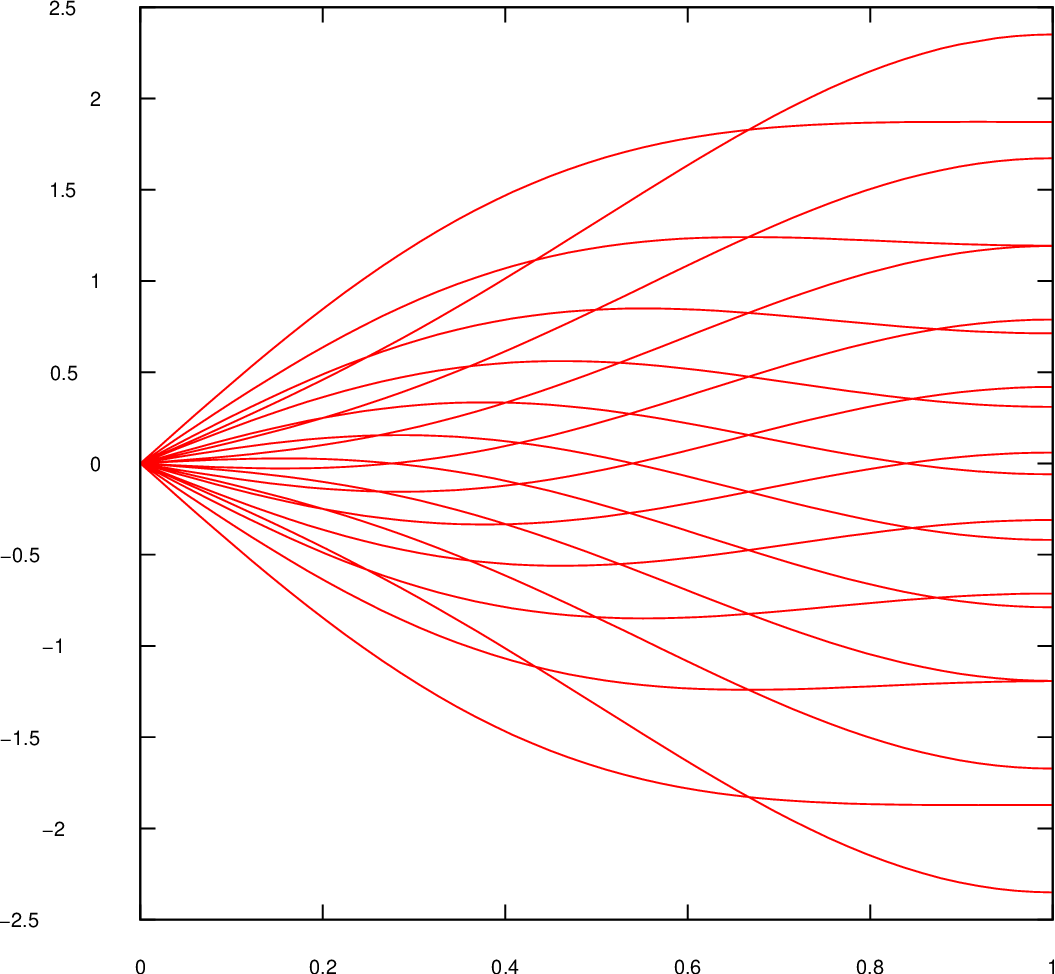}
\end{minipage} \hfill
\begin{minipage}[c]{.46\linewidth}
\psfrag{-0.2}{$-0.2$}
\psfrag{-0.4}{$-0.4$}
\psfrag{-0.6}{$-0.6$}
\psfrag{-0.8}{$-0.8$}
\psfrag{-1}{$-1$}
\psfrag{0}{$0$}
\psfrag{0.2}{$0.2$}
\psfrag{0.4}{$0.4$}
\psfrag{0.6}{$0.6$}
\psfrag{0.8}{$0.8$}
\psfrag{1}{$1$}
\includegraphics[height=5.5cm]{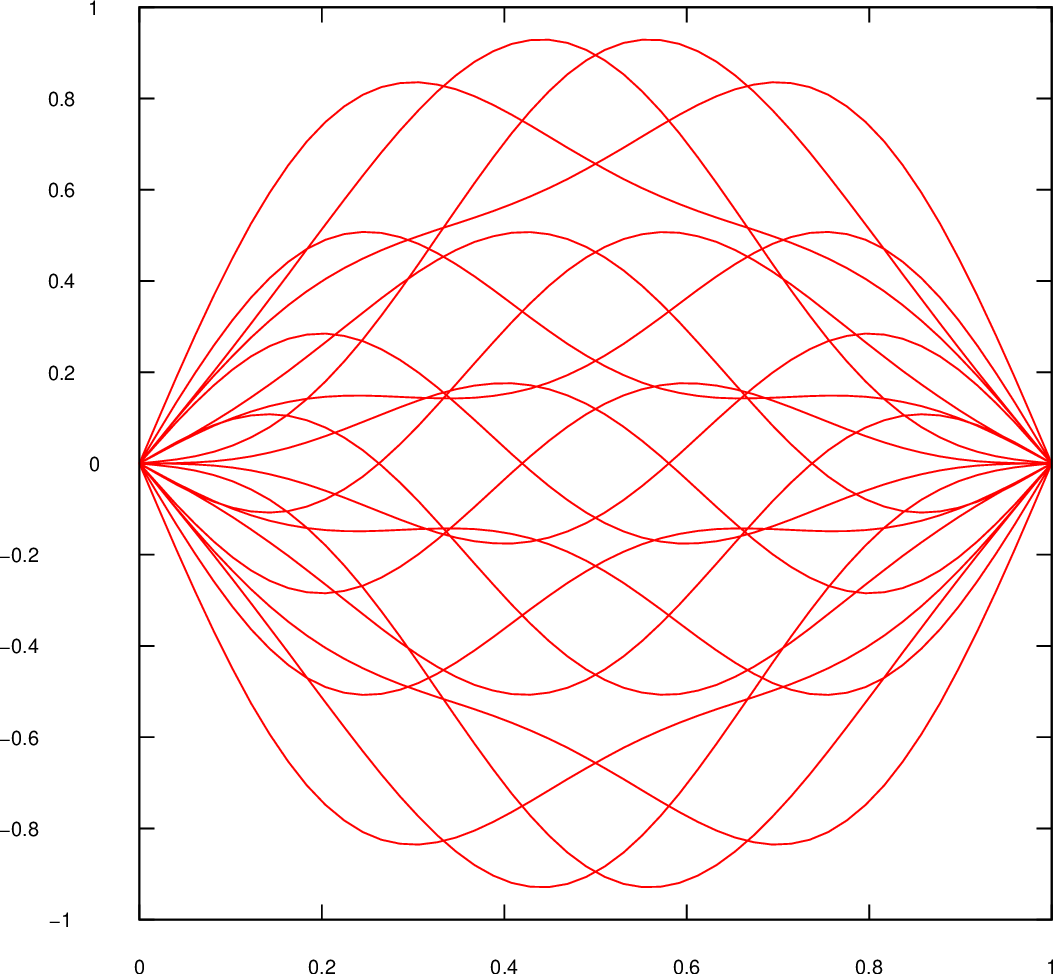}
\end{minipage}
\caption{Optimal product quantization of Brownian motion (left) and Brownian bridge (right) on $[0,1]$.}
\label{fig:brownian_motion_product_quantizer}
\end{figure}

\begin{figure}[!ht]
\begin{minipage}[c]{.46\linewidth}
	\psfrag{-1.5}{$-1.5$}
	\psfrag{-1}{$-1$}
	\psfrag{-0.5}{$-0.5$}
	\psfrag{0}{$0$}
	\psfrag{0.5}{$0.5$}
	\psfrag{1}{$1$}
	\psfrag{1.5}{$1.5$}
	\psfrag{2}{$2$}
	\psfrag{2.5}{$2.5$}
	\psfrag{3}{$3$}
\includegraphics[height=5.5cm]{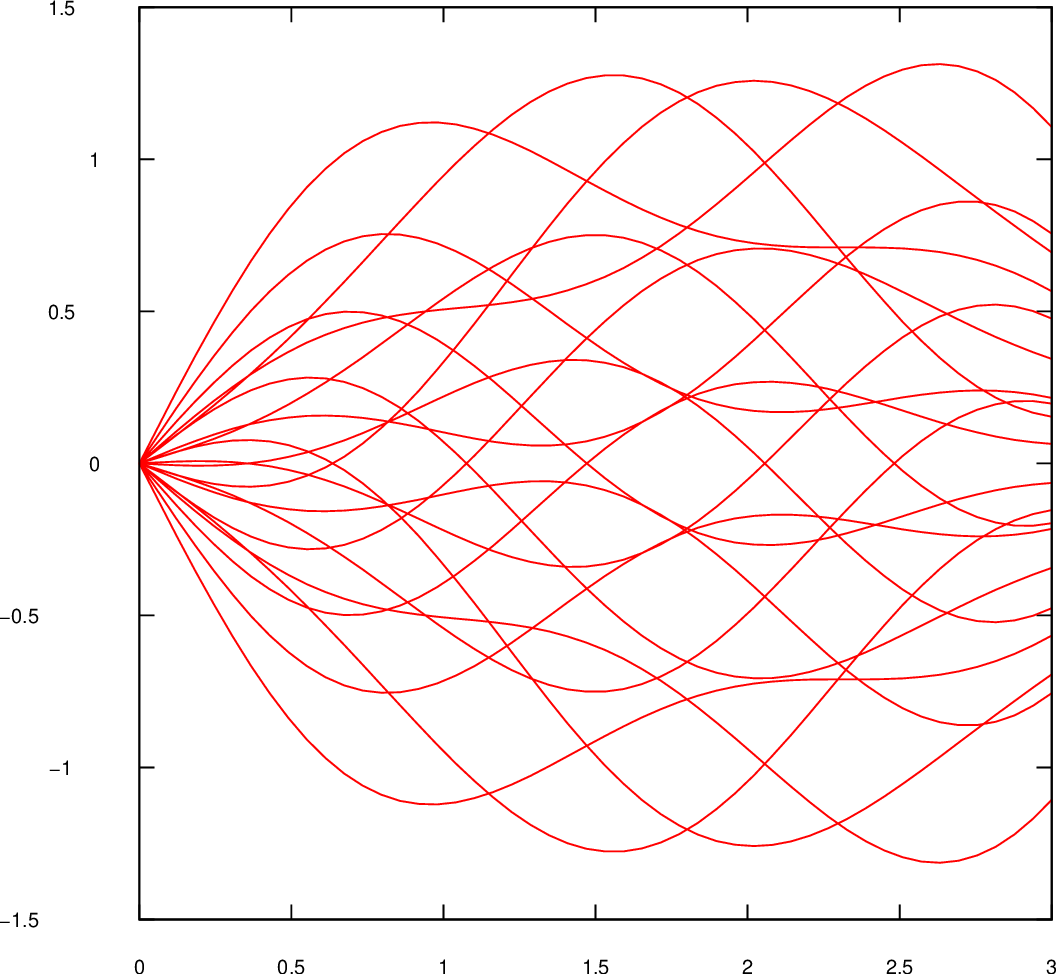}
\end{minipage} \hfill
\begin{minipage}[c]{.46\linewidth}
	\psfrag{-1.5}{$-1.5$}
	\psfrag{-1}{$-1$}
	\psfrag{-0.5}{$-0.5$}
	\psfrag{0}{$0$}
	\psfrag{0.5}{$0.5$}
	\psfrag{1}{$1$}
	\psfrag{1.5}{$1.5$}
	\psfrag{2}{$2$}
	\psfrag{2.5}{$2.5$}
	\psfrag{3}{$3$}
\includegraphics[height=5.5cm]{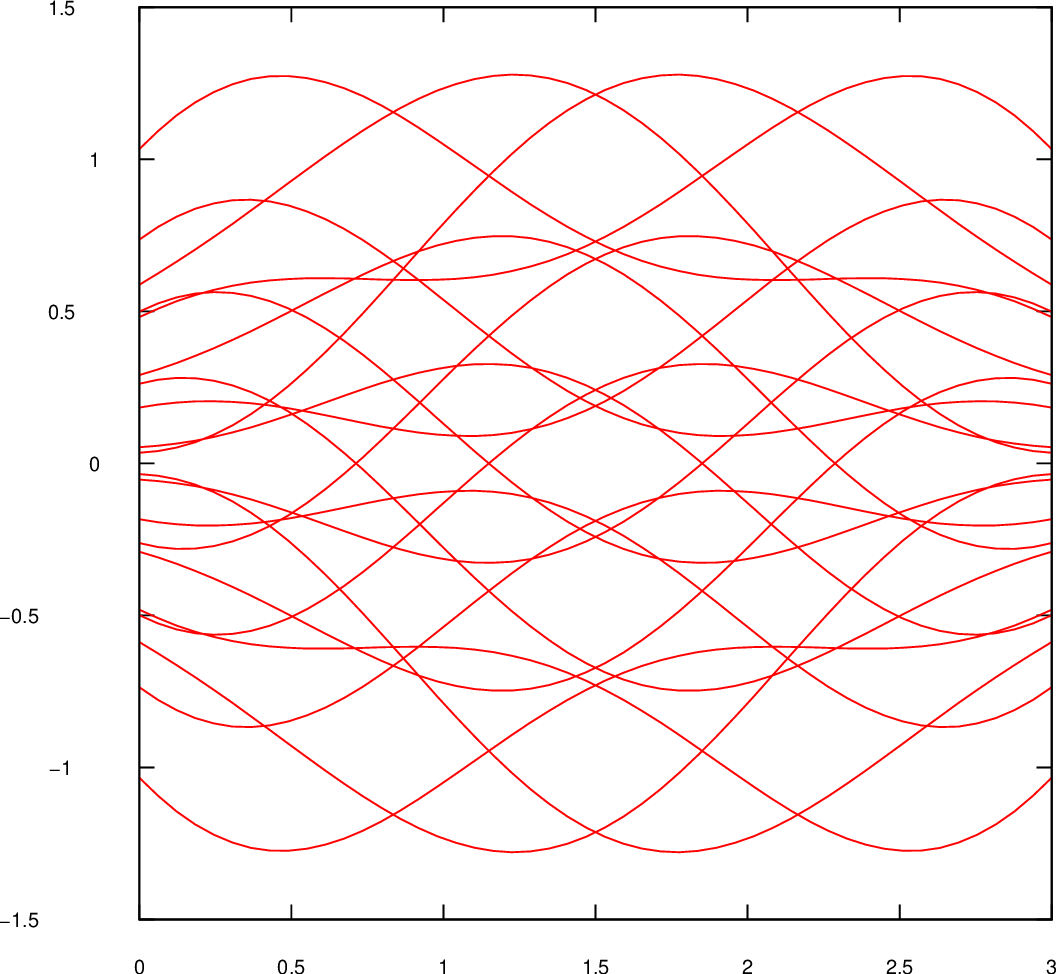}
\end{minipage}
\caption{Optimal product quantization of a centered Ornstein-Uhlenbeck process, starting from $X_0 = 0$ (left) and stationary (right) solution of the SDE $d X_t = -X_t dt + dW_t$, on $[0,3]$.}
\label{fig:ornstein_uhlenbeck_product_quantizer}
\end{figure}

 In Figure \ref{fig:brownian_motion_product_quantizer}, we display optimal product quantizers of Brownian motion and Brownian bridge on $[0,1]$. In Figure \ref{fig:ornstein_uhlenbeck_product_quantizer}, we display optimal product quantizers of the centered Ornstein-Uhlenbeck process starting from $X_0=0$ and a stationary Ornstein-Uhlenbeck on $[0,3]$.
\subsubsection{Rate of decay of the quantization error}
\par The rate of decay of the quadratic functional quantization error of Gaussian processes was first investigated in \cite{LuschgyPagesFunctional3} and more precise results were then established in \cite{LuschgyPagesFunctional2}. These results rely on assumptions on the asymptotic behavior of the Karhunen-Loève eigenvalues of the considered process. 
\par \noindent Let $X$ be a bi-measurable centered Gaussian process on $[0,T]$ of continuous covariance function $\Gamma^X$ and such that $\int_0^T \E[X_s^2]ds <\infty$. Its Karhunen-Loève eigensystem is denoted by $\left(e^X_n,\lambda^X_n\right)_{n \geq 1}$.
\begin{theo}[Quadratic quantization error asymptotics]\label{thm:KL_convergence}
\par Assume that $\lambda^X_n \sim \phi(n)$ as $n \to \infty$, where $\phi : (s,\infty)\to (0,\infty)$ is a decreasing function such that $\lim\limits_{x \to \infty} \frac{\phi(tx)}{\phi(x)} = t^{-b}$ for $b>1$ and $s>0$. Set $\psi(x) := \frac{1}{x \phi(x)}$. Then
$$
\mathcal{E}_N(X) \sim \left( \left( \frac{b}{2}\right)^{b-1} \frac{b}{b-1} \right)^{1/2} \psi(\log(N))^{-1/2} \quad \textnormal{as } N \to \infty.
$$
Moreover, the optimal product quantization dimension $m^X(N)$ verifies $m^X(N) \sim \frac{2}{b} \log(N) \quad \textnormal{as } N \to \infty$, and the optimal product quantization error $\mathcal{E}^{pq}_N(X)$ of level $N$ satisfies
$$
\mathcal{E}^{pq}_N(X) \lesssim \left( \left( \frac{b}{2}\right)^{b-1} \frac{b}{b-1} + C(1) \right)^{1/2} \psi(\log(N))^{-1/2} \quad \textnormal{as } N \to \infty,
$$
where $C(1)$ is a universal positive constant. 
\end{theo}
\par \noindent A proof is available in \cite{LuschgyPagesFunctional2}. Despite of the fact that optimal product quantization is not asymptotically optimal, it provides a rate-optimal sequences of quantizers. Typical rates are $\underset{N \to \infty}{\sim} \log(N)^{-\alpha}$ for $\alpha>0$. For Brownian motion, Brownian bridge and Ornstein-Uhlenbeck processes, we have $\alpha = \frac{1}{2}$. 
\section[Quantization as a control variate]{A first attempt to quantization-based variance reduction: quantization as a control variate}\label{seq:quantization_variance_reduction}
\par This method has been originally proposed in \cite{PagesPrintemsFunctional4}. Let $X$ be an $E$-value $L^2$ random variable, consider $N \in \N^*$ and let $\Gamma = \{y_1, \ldots, y_N\}$ be an $N$-codebook. We define a quantizer $Y$ of $E$ by $Y := \Proj(X)= \sum\limits_{i=1}^N y_i \1_{C_i}(X)$ where $C=\{C_1, \ldots,C_N\}$ is a partition of $E$. At this stage, we do not need $\Proj$ to be a nearest neighbor projection onto $\Gamma$. 
\par Let $F : E \to E$ be a Lipschitz continuous function. In order to compute $\E[F(X)]$, we use that:
\begin{equation}\label{eq:control_variate_quantization}
	\begin{array}{lll}
	\E[F(X)] 	&= \E\left[F(\Proj(X))\right] + \E\left[F(X)-F(\Proj(X))\right] \\
				&= \underbrace{\E\left[F(\Proj(X))\right]}_{(a)} + \underbrace{\frac{1}{M} \sum\limits_{m=1}^M F\left(X^{(m)}\right)-F\left(\Proj\left(X^{(m)}\right)\right)}_{(b)} + R_{N,M},
	\end{array}
\end{equation}
where $X^{(m)}, 1 \leq m \leq M$ are $M$ independent copies of $X$, and $R_{N,M}$ is a remainder term defined by Equation \eqref{eq:control_variate_quantization}. Term $(a)$ is computed by quantization-based cubature and Term $(b)$ is computed by a Monte Carlo simulation. We have 
$$
\|R_{N,M}\|_2 = \frac{\sigma(F(X)-F(\Proj(X)))}{\sqrt{M}} \leq \frac{\|F(X)-F(\Proj(X)) \|_2}{\sqrt{M}}  \leq [F]_{\textnormal{Lip}} \frac{\|X-\Proj(X)\|_2}{\sqrt{M}}.
$$
Furthermore, $\sqrt{M} R_{N,M} \ \overset{\mathcal{L}}{\rightarrow}\ \mathcal{N}\Big(0,\var\big(F(X)-F(\Proj(X))\big)\Big)$.
\par \noindent \textbf{Consequently, in the $d$-dimensional case}, if $F$ is Lipschitz continuous and $\left(\widehat{X}^N\right)_{N\in \N} = (\Proj^N(X))_{N \in \N}$ is a rate-optimal sequence of quantizers of $X$, then we have $\left\|F(X)-F\left(\Proj^N(X)\right)\right\|_2 \leq [F]_{\textnormal{Lip}} \frac{C_X}{N^{1/d}}$ so that
$$
\|R_{N,M}\|_2 \leq [F]_{\textnormal{Lip}} \frac{C_X}{M^{1/2} N^{1/d}}. 
$$
\par \noindent \textbf{Likewise, in the case of Brownian motion}, if $\left(\widehat{W}^N\right)_{N\geq 1}$ is a rate-optimal sequence of quadratic K-L product quantizers of Brownian motion, if $F$ is a Lipschitz continuous functional, then
$\left\|F(W) - F\left(\widehat{W}^N\right)\right\|_2 \leq [F]_{\textnormal{Lip}} \frac{C_W}{\log(N)^{1/2}}$
so that
$$
\left\|R_{N,M}\right\|_2 \leq [F]_{\textnormal{Lip}} \frac{C_W}{M\log(N)^{1/2}}.
$$
\subsubsection*{The bottleneck of fast nearest neighbor search}
\par $\bullet$ \textbf{The complexity of the projection:} When implementing the quantization-based control variate variable method \eqref{eq:control_variate_quantization} , for every draw of the Monte Carlo simulation, one has to compute the projection $\Proj(X^{(m)})$. As a consequence, the efficiency of the method is conditioned by the efficiency of the projection procedure. When dealing with Voronoi quantization, this is simply the nearest neighbor projection. 
\par The problem of nearest neighbor projection, also known as the post-office problem \cite{KnuthArt3}, has been widely investigated in the area of computational geometry. It has been solved near optimally in the low dimensional case. Algorithms differ on their practical efficiency on real data sets. For large dimensions, most solutions have a complexity that is exponential with the dimension, or require a longer query time than the obvious brute force algorithm. In fact for dimension $d>\log N$, a brute force algorithm is usually the best choice. Still, even in low dimension, fast nearest neighbor search is a critical part of the algorithm. Let us mention \cite{CorlayFNNS} for a fast nearest neighbor search algorithm based on recursive vector quantization. 
\par The speed of the projection can also be increased by relaxing the hypothesis that the projection onto the quantizer is a nearest neighbor projection or by choosing simpler partitions of the state space. 
\vspace{1mm}\\
$\bullet$ \textbf{The functional case:} The problem of nearest neighbor search is even less tractable in the functional case, as one does not simulate the whole trajectory of the stochastic process but only its marginals at discrete dates, and therefore we can only make an assumption on the interpolation to compute the nearest neighbor.
\par In \cite{LejayReutenauerControl}, a functional quantizer of Brownian motion is used as a control variate variable.
\section{Application of quantization to stratification}\label{sec:stratification}
\subsection{Some background on stratified sampling}
\par The main idea of stratification is to localize the Monte Carlo simulation on the elements of a measurable partition of the state space of an $L^2$ random variable $X:(\Omega,\mathcal{A}) \to (E,\mathcal{E})$. Let $(A_i)_{i \in I}$ be a finite $\mathcal{E}$-measurable partition of $E$. The sets $A_i$ are called \emph{strata}. We assume that the weights $p_i = \PP\left(X \in A_i \right)$, $i \in I$ are positive. We will make two pseudo or operating assumptions on these strata:
\begin{itemize}
\item $\forall i \in I$, $p_i = \PP(X \in A_i)$ is known. 
\item $\forall i \in I$, the random variable $X_i \simdist \mathcal{L}(X |X \in A_i)$ can be simulated at a reasonable cost (say similar to that of $X$ itself).
\end{itemize}
\par Tractability of simulation is a major constraint for practical implementation and it has a strong impact on the design of the strata. In practice, we can formulate the condition by assuming that $X_i= \phi_i(U)$ where $U$ is uniformly distributed on $[0,1]^{r_i}$ and $\phi_i:[0,1]^{r_i} \to \R$ is an easily computable function. (We have $r_i \in \N \cup \{+\infty\}$, the case $r_i = + \infty$ corresponds to the acceptance-rejection method.)
\vspace{3mm}
\par Let $F:(E,\mathcal{E}) \to (\R,\mathcal{B}(\R))$ such that $\E[|F(X)|]<+\infty$. We have
$$
\E[F(X)] = \sum\limits_{i \in I} \E[\1_{\{ X \in A_i\}}F(X)] = \sum\limits_{i \in I} p_i \E[F(X) | X \in A_i] = \sum\limits_{i \in I} p_i \E[F(X_i)].
$$
\par The stratification concept comes into play now. Let $M$ be the global budget allocated to the computation of $\E[F(X)]$ and let $M_i=q_i M$ be the budget allocated to compute $\E[F(X_i)]$ in each stratum (with $0 \leq q_i \leq 1$, $i \in I$ and $\sum\limits_{i \in I} q_i = 1$). This leads to define the (unbiased) estimator of $\E[F(X)]$:
\begin{equation}\label{eq:unbiased_stratif_estimator}
\overline{F(X)}^I_M := \sum\limits_{i \in I} p_i \frac{1}{M_i} \sum\limits_{k=1}^{M_i} F\left(X_i^k\right),
\end{equation}
where $(X_i^k)_{1\leq k \leq M_i}$ is a $\mathcal{L}(X|X \in A_i)$-distributed random sample. We have
\begin{equation}\label{eq:stratif_estimator_variance}
\var\left(\overline{F(X)}^I_M\right)=\frac{1}{M} \sum\limits_{i\in I} \frac{p_i^2}{q_i} \sigma_{F,i}^2,
\end{equation}
where $\sigma^2_{F,i} = \var(F(X)|X\in A_i) = \var(F(X_i))$, $i \in I$. Optimizing the allocation of the number of draws to the different strata amounts to solving the following minimization problem: 
\begin{equation}\label{eq:stratification_minimization}
\min\limits_{(q_i) \in \mathcal{P}_I} \sum\limits_{i \in I} \frac{p_i^2}{q_i} \sigma^2_{F,i} \hspace{2mm} \textnormal{ where } \mathcal{P}_I := \left\{ (q_i)_{i\in I} \in \R_+^I \middle| \sum\limits_{i\in I} q_i = 1 \right\}. 
\end{equation}
\subsubsection{Natural stratified sampling}\label{sec:sub_optimal_choice}
A natural choice is to set
\begin{equation}\label{eq:sub_optimal_choice}
q_i=p_i, \hspace{5mm} i \in I.
\end{equation}
since the weights $p_i$ are known. Furthermore, this always reduces the variance. 
\begin{multline*}
\sum\limits_{i\in I} \frac{p_i^2}{q_i} \sigma^2_{F,i} = \sum\limits_{i\in I} p_i \sigma^2_{F,i} = \sum\limits_{i\in I} \E\left[\Big( F(X)-\E[F(X)|X \in A_i]\Big)^2 \1_{A_i}(X) \right]\\
= \|F(X) - \E[F(X)|\sigma(\{X \in A_i \}, \ i\in I)] \|_2^2\\
\leq \|F(X) - \E[F(X)] \|_2^2 = \var(F(X)).
\end{multline*}
\subsubsection{Optimal stratified sampling}\label{sec:optimal_choice}
\par The optimal choice is the solution to the constrained minimization problem \eqref{eq:stratification_minimization}. Schwarz's inequality yields
$$
\sum\limits_{i \in I} p_i \sigma_{F,i} = \sum\limits_{i \in I} \frac{p_i \sigma_{F,i}}{\sqrt{q_i}} \sqrt{q_i} \leq \Big(\sum\limits_{i \in I} \frac{p_i^2 \sigma_{F,i}^2}{q_i} \Big)^{1/2}\Big( \sum\limits_{i\in I} q_i \Big)^{1/2}.
$$
\par \noindent The solution corresponds to the equality case in Schwarz's inequality, that is 
\begin{equation}\label{eq:optimal_budget}
q_i^* = \frac{p_i \sigma_{F,i}}{\sum\limits_{j\in I} p_j \sigma_{F,j}}, \quad i\in I
\end{equation}
with a resulting minimal variance of $\Big(\sum\limits_{i \in I} p_i \sigma_{F,i}\Big)^2$. At this stage, the problem is that we do not \textit{a priori} know the local inertia $\sigma_{F,i}^2$. Still, using the fact that $L^p$ norms are decreasing with $p$, we see that 
$$
\sigma_{F,i} \geq \E \Big[ \left|F(X)-\E\left[F(X)\middle|\{X \in A_i \}\right] \right| \Big| \{ X \in A_i \}\Big],
$$
so that
$$
\left(\sum\limits_{i \in I} p_i \sigma_{F,i}\right)^2 \geq \Big\|F(X)-\E\left[F(X)\middle|\sigma(\{X \in A_i\}, \ i\in I)\right] \Big\|_1^2.
$$
\par In \cite{JourdainStratification1}, Étoré and Jourdain proposed an algorithm which adaptively modifies the proportion of further drawings in each stratum and which converges to the optimal allocation. 
\par In Section \ref{sec:quantization_based_stratif}, we show that the problem of designing good strata, in term of variance reduction is linked with optimal quantization. Besides, with quantization-based stratified sampling, the weights $p_i$ are already known.
\subsection{Quantization and stratified sampling}\label{sec:quantization_based_stratif}
\par The main drawback of using quantization as a control variate is the repeated computations of the projections onto the quantizer. (Nearest neighbor searches in the case of a Voronoi quantizer.) In the case of stratified sampling, \emph{one does not have to use a projection procedure}. Instead, we must focus on the cost of the simulation of conditional distributions $\mathcal{L}(X|X\in A_i)$, $i \in I$.
\par Proposition \ref{prop:universal_stratification} brings together previous results and highlights the relationships with quantization. It shows that stratification has uniform efficiency among the class of Lipschitz continuous functionals. 
\begin{prop}[Universal stratification]\label{prop:universal_stratification}
\par Let $A=(A_i)_{i\in I}$ be a partition of $E$ and let $\Proj_{A,X}$ denote the centroidal projection associated with $X$ and $A$, defined in Definition \ref{def:centroidal_projection}. 
\begin{enumerate}
\item\label{it:us1} Considering the local inertia of $X$ in $A_i$, $\sigma_i^2 = \E\left[\left|X-\E[X|X\in A_i]\right|^2 \middle| X \in A_i \right]$, we have for every Lipschitz continuous function $F:E \to E$, $\sigma_{F,i} \leq [F]_{\textnormal{Lip}} \sigma_i$ where $[F]_{\textnormal{Lip}} = \sup\limits_{x \neq y} \frac{F(x)-F(y)}{|x-y|}$, so that
	\begin{equation}\label{eq:strat_ineq1}
	\sup\limits_{[F]_{\textnormal{Lip}} \leq 1} \sigma_{F,i} = \sigma_i,
	\end{equation}
\item\label{it:us2}	In the case of natural stratified sampling (see Section \ref{sec:sub_optimal_choice}),
		\begin{equation}\label{eq:strat_ineq2}
		\sup\limits_{[F]_{\textnormal{Lip}} \leq 1} \Big( \sum\limits_{i \in I} p_i \sigma_{F,i}^2 \Big) = \sum\limits_{i\in I} p_i \sigma_i^2 = \Big\|X-\E[X|\sigma(\{ X \in A_i \}, \ i\in I)] \Big\|_2^2 = \Big\|X-\Proj_{A,X}(X) \Big\|_2^2.
		\end{equation}
\item\label{it:us3}	In the case of the optimal choice (see Section \ref{sec:optimal_choice}), 
		\begin{equation}\label{eq:strat_ineq3}
		\sup\limits_{[F]_{\textnormal{Lip}} \leq 1} \Big( \sum\limits_{i \in I} p_i \sigma_{F,i} \Big)^2 = \Big( \sum\limits_{i\in I} p_i \sigma_i \Big)^2,
		\end{equation}
		and
		$$
		\Big( \sum\limits_{i\in I} p_i \sigma_i \Big)^2 \geq \Big\|X-\E[X|\sigma(\{ X \in A_i \}, \ i\in I)] \Big\|_1^2 = \Big\|X-\Proj_{A,X}(X) \Big\|_1^2.
		$$
\item\label{it:us4} In the case of real-valued Lipschitz continuous functions $F : E \to \R$, Equalities \eqref{eq:strat_ineq1}, \eqref{eq:strat_ineq2} and \eqref{eq:strat_ineq3} hold as inequalities. 
\end{enumerate}
\end{prop}
\noindent \textbf{Proof:} We have 
\begin{multline*}
\sigma_{F,i}^2 = \var\left(F(X)\middle|X\in A_i\right) = \E\left[ \left|F(X)-\E[F(X)|X\in A_i]\right|^2 \middle| X \in A_i \right]\\
\leq \E\left[ \left|F(X)-F(\E[X|X\in A_i])\right|^2 \middle| X \in A_i \right].
\end{multline*}
Now using that $F$ is Lipschitz continuous, we get 
$$
\sigma_{F,i}^2 \leq [F]^2_{\textnormal{Lip}} \frac{1}{p_i} \E \left[ \left|X-\E[X|X\in A_i]\right|^2 \1_{\{X\in A_i\}} \right] = [F^2]_{\textnormal{Lip}} \sigma_i^2.
$$
\par \noindent Items \ref{it:us2} and \ref{it:us3} easily follow from Item \ref{it:us1}. Equality follows by considering $F = Id_E$. \myqed
\subsubsection{Universal stratified sampling}
\par \noindent Proposition \ref{prop:universal_stratification} suggests, in the case of Lipschitz continuous functionals, to set 
$$
q_i = \frac{p_i \sigma_i}{\sum\limits_{j \in I} p_j \sigma_j}, \quad j \in I,
$$
so that we have uniform efficiency among the class of Lipschitz continuous functionals. This allocation scheme will be further referred to as the ``universal stratification'' weights. It also shows that, in the Lipschitz continuous case, it is always beneficial to reduce the quadratic distortion associated with the centroidal projection $\Proj_{A,X}$. 

\par Still, this minimization should not be done at the expense of the efficiency of the simulation of the corresponding conditional distributions. We should reach for a balance between the efficiency of the simulation in the strata and the quadratic quantization error controlling the variance reduction. For example, in Section \ref{sec:functional_stratification}, in the functional case, we will use optimal product quantizers, which are rate optimal (and numerically near optimal) and allow for a much more efficient simulation than real optimal functional quantization.
\begin{remark} 
\par \noindent We should also mention the adaptive strata design proposed in \cite{JourdainAdaptiveStrat, JourdainConvenientStratDir}.
\end{remark}
\subsection{Simulation in hyper-rectangular strata in the independent Gaussian case}\label{seq:gaussian_rectangle}
\par Consider $X \simdist \mathcal{N}(0,I_d)$, $d\geq 1$ and $(e_1, \ldots, e_d)$ an orthonormal basis of $\R^d$. Let $N_1, \ldots, N_d \geq 1$ be the number of strata in each direction and for $1\leq i \leq d$, $-\infty = \gamma_0^i \leq \gamma_1^i \leq \cdots \leq \gamma_{N_i}^i = + \infty$. We define
$$
A_{\ui} := \bigcap\limits_{l=1}^d \Big\{ x \in \R^d \textnormal{ such that } \langle e_l,x \rangle \in [x^l_{i_l-1},x^l_{i_l}] \Big\}, \hspace{5mm }\ui \in \prod\limits_{l=1}^d \{1, \ldots, N_l \}.
$$
Then for every $\ui \in \prod\limits_{l=1}^d \{1, \ldots, N_l \}$, 
$\mathcal{L}\left(X \middle| X \in A_{\ui} \right) = \bigotimes_{l=1}^d \mathcal{L}\left(Z \middle| Z \in \left[\gamma^l_{i_l-1},\gamma^l_{i_l}\right]\right), \textnormal{ where }Z \simdist \mathcal{N}(0,1)$, $p_{\ui} = \PP(A_{\ui}) = \prod\limits_{l=1}^d \big( \mathcal{N}(\gamma_{i_l}^l) - \mathcal{N}(\gamma_{i_{l}-1}^l) \big)$ and for $-\infty \leq a \leq b \leq \infty$, 
\begin{equation}\label{eq:gauss_trunc_closed}
\mathcal{L}\left(Z \middle| Z \in [a,b]\right) = \mathcal{N}^{-1} \left( \left(\mathcal{N}(b) - \mathcal{N}(a)\right) U + \mathcal{N}(a) \right), \hspace{5mm} U \simdist \mathcal{U}([0,1]).
\end{equation}
\section{Functional stratification of Gaussian processes}\label{sec:functional_stratification}
\par In the functional case, the state space of the random values are functional spaces. What is usually done is to simulate a scheme to approximate marginals of the underlying process. 
\par In this section, we assume that $X$ is a centered $\R$-valued bi-measurable Gaussian process on $[0,T]$ that satisfies $\int_0^T \E[X_t^2]dt < \infty$. We are interested by the value of 
$\E[F(X_{t_0}, X_{t_1}, \ldots, X_{t_n})]$ for some real function $F$, where $0 = t_0 \leq t_1 \leq \cdots \leq t_n = T$ are $n+1$ dates of interest for the underlying process. 
\par (For example, $X$ can be a standard Brownian motion on $[0,T]$, and one computes the risk-neutral expectation of a path-dependent payoff of a diffusion based on $X$.)
\par The results of this section can be easily generalized to the multi-dimensional case, like multifactor diffusions. Still we restrict ourselves to the one-dimensional setting for clarity. 
\par Let us assume that $\chi \in \mathcal{O}_{pq}(X,N)$ is a K-L optimal product quantizer of $X$. The codebook associated with this product quantizer is the set of the paths of the form 
$$
\chi_{\ui} = \sum\limits_{n \geq 1} \sqrt{\lambda_n^X} \gamma_{i_n}^{(N_n)} e_n^X, \hspace{5mm} \ui = \{i_1, \ldots, i_n, \ldots \},
$$
with the same notation as in Section \ref{sec:functional_product_quantization}. We now need to be able to simulate the conditional distribution 
$$
\mathcal{L}(X | X \in A_{\ui})
$$
where $A_{\ui}$ is the cell associated with $\chi_{\ui}$ in the codebook. To simulate the conditional distribution $\mathcal{L}(X | X \in A_{\ui})$, one will :
\begin{itemize}
	\item First, simulate the first K-L coordinates of $X$, using \eqref{eq:gauss_trunc_closed}.
	\item Then simulate the conditional distribution of the marginals of the Gaussian process given its first K-L coordinates. 
\end{itemize}
\begin{remark}
\par We have chosen to use K-L optimal product quantizers instead of optimal quantizers because in this case, the Voronoi cells in this are hyper-rectangles, which allows us to simulate the first K-L coordinates more easily than in the general case. Moreover, the rate of decay of the quantization errors is rate-optimal under some conditions on the Karhunen-Loève eigenvalues which are verified in the considered examples \cite{LuschgyPagesFunctional3}.
\end{remark}
\subsection{Simulation of marginals of the Gaussian process, given its $d$ first K-L coordinates}\label{sec:marginal_simulation}
\par In this setting, the aim is to simulate the conditional distribution 
\begin{equation}\label{eq:conditional_distribution}
\mathcal{L}\Big(X_{t_0}, \ldots, X_{t_n} \Big| \int_0^T X_s e_1^X(s) ds, \int_0^T X_s e_2^X(s) ds, \ldots , \int_0^T X_s e_d^X(s) ds \Big)
\end{equation}
where $(X_t)_{t \in [0,T]}$ is an $L^2$ $\R$-valued Gaussian process, and $(e_k^X,\lambda_k^X)_{k\in \N^*}$ is the Karhunen-Loève system associated with the process $X$. Hence $\Big( X_{t_0}, \ldots, X_{t_n}, \int_0^T X_s e_1^X(s)ds, \ldots , \int_0^T X_s e_d^X(s) ds \Big)$ is a Gaussian vector. As a consequence, if we denote $Y:=\left(\begin{array}{ccc} \int_0^T X_s e_1^X(s) ds \\ \vdots \\ \int_0^T X_s e_d^X(s) ds \end{array}\right)$ and $V:=\left(\begin{array}{ccc} X_{t_0} \\ \vdots \\ X_{t_n} \end{array}\right)$, the conditional distribution (\ref{eq:conditional_distribution}) is given by the transition kernel 
$\nu(y,A) = \mathcal{N}\left(Af_{V|Y}(y),\cov(V-\E[V|Y])\right)$, where $Af_{V|Y} : \R^d \to \R^n$ is an affine function corresponding to the linear regression of $V$ on $Y$, $Af_{V|Y}(Y) := \E[V|Y]$. 
\begin{itemize}
\item We have $Af_{V|Y}(Y) = \E[V] + R_{V|Y} Y$ where $R_{V|Y} = \cov(V,Y) \cov(Y)^{-1}$. Using that $\cov(Y) = \Big(\lambda_i^X \delta_{ij}\Big)_{1\leq i,j \leq d}$ and $\cov(V,Y) = (\lambda_k^X e_k^X(t_i))_{0\leq i \leq n, 1 \leq k \leq d}$, we get
\begin{equation}\label{eq:r1_value}
R_{V|Y} = \left(e_j^X(t_i)\right)_{0\leq i \leq n, 1 \leq j \leq d}.
\end{equation}
\item The covariance matrix is 
\begin{multline*}
K:= \cov\left(V-\E[V|Y]\right)	= \E\left[\left(V-R_{V|Y}Y\right)\left(V-R_{V|Y}Y\right)\right]\\
		= \cov(V) - 2 \cov\left(V,R_{V|Y}Y\right) + \cov\left(R_{V|Y} Y\right) = \cov(V) - \cov\left(R_{V|Y} Y\right)\\
				= \left( \cov(V_l,V_k) - \sum\limits_{i=1}^d \lambda_i e_i^X(t_l) e_i^X(t_k) \right)_{0 \leq k,l \leq n}.
\end{multline*}
\end{itemize}
\par The easiest way to simulate according to this probability distribution would be to use the Cholesky factorization of $K$. However, when using this method, the simulation of a simple path involves the quadratic complexity of an $n \times n$ matrix multiplication, which is not satisfactory for our purpose. 
\subsection{Faster simulation of conditional paths - Bayesian simulation}\label{sec:faster_simulation}
\par As pointed out earlier, the naive simulation method for $\mathcal{L}(V|Y)$ requires for each path a multiplication by a Cholesky transform of $K$ whose cost is $O(n^2)$.  
\begin{itemize}
\item Yet, the quantization dimension $d$ of the process is close to $\log(N)$ where $N$ is the number of strata, and $n$, the number of time steps, is usually very large compared to $d$. 
\item The idea here is that the conditional distribution $\mathcal{L}(V|Y)$ is determined through the Bayes lemma, by the conditional distribution $\mathcal{L}(Y|V)$ and the two marginal distributions $\mathcal{L}(V)$ and $\mathcal{L}(Y)$. 
\end{itemize}
\par One knows that $V = \E[V|Y] \ {\overset{\independent }{+}} \ Z$ where $Z \simdist \mathcal{N}(0,\cov(V-\E[V|Y]))$ is independent of $Y$. Hence one is able to simulate according to $\mathcal{L}(V|Y=y)$ if one can simulate the distribution of $Z$, writing $\mathcal{L}(V|Y=y) = \E[V|Y=y] + \mathcal{L}(Z)$. This decomposition corresponds to the splitting of the Karhunen-Loève expansion:
$$
\left(\begin{array}{c}
V_0\\
\vdots\\
V_n
\end{array}\right)
= \underbrace{\sum\limits_{k=1}^d \underbrace{\sqrt{\lambda_k^X} \xi_k}_{= Y_k} 
\left(\begin{array}{c}
e_k^X(t_0)\\
\vdots\\
e_k^X(t_n)
\end{array}\right)}_{=\E[V|Y]}
\ {\overset{\independent }{+}} \
\underbrace{\sum\limits_{l \geq d+1}
\sqrt{\lambda_k^X} \xi_k
\left(\begin{array}{c}
e_k^X(t_0)\\
\vdots\\
e_k^X(t_n)
\end{array}\right)}_{=Z}.
$$
\par \noindent To simulate $Z$, one simulates the distribution of $V$ and the conditional distribution $\mathcal{L}(Z|V)$. 
$$
\begin{array}{lll}
\textnormal{We have} \hspace{7mm} \mathcal{L}(Z|V)	& \simdist \delta_V- \mathcal{L}(\E[V|Y]|V) \simdist \delta_V- Af_{V|Y} \mathcal{L}(Y|V) \\
												& \simdist \delta_V- Af_{V|Y} \mathcal{N}(\E[Y|V],\cov(Y-\E[Y|V])). \\
\end{array}
$$
 	If $Af_{Y|V}$ is the affine function corresponding to the regression of $Y$ on $V$ and $R_{Y|V}$ its linear part,
	$$
	\cov(Y-\E[Y|V]) = \cov(Y) + \cov(\E[Y|V]) -2 \cov(Y,\E[Y|V]) = \cov(Y) - R_{Y|V} \cov(V) ^t\!R_{Y|V}.
	$$
This yields $Z = V - Af_{V|Y} (G)$ where $G \simdist \mathcal{N}(Af_{Y|V}(V),\cov(Y) - R_{Y|V} \cov(V) ^t\!R_{Y|V})$. Finally, we can use the following method to simulate the conditional distribution of $V$  $Y$. 
\begin{breakbox}
	\begin{itemize}
		\item Simulate $V$. \hfill \emph{$O(n)$.}
		\item Simulate $G \simdist \mathcal{N}\left(Af_{Y|V}(V),\cov(Y) - R_{Y|V} \cov(V) ^t\!R_{Y|V}\right)$ \hfill \emph{$O(d \times d)$}.
		\item Compute $Z = V - Af_{V|Y} (G)$. \hfill \emph{$O(d \times n)$}.
	\end{itemize}
	\begin{center}
	The random variable $T := Af_{V|Y}(y) + Z$ satisfies $T \simdist \mathcal{L}(V|Y=y)$.
	\end{center}
\end{breakbox}
\vspace{2mm}
\par Let us remind that $Af_{V|Y}$ is trivially defined in Equation (\ref{eq:r1_value}), because coordinates of $Y$ are independent. Other matrices implied in this algorithm are computed prior to any Monte Carlo simulation. In general, $R_{Y|V}$ can simply be computed by performing a numerical least-square regression. Moreover in the special cases of Brownian motion, Bronian bridge and Ornstein-Uhlenbeck processes, there are closed-form expressions for the $R_{Y|V}$, which we present in Appendix \ref{sec:annex_R2_computation}. 
\par \noindent In the case of Brownian motion, for a uniform time discretization mesh $t_j = \frac{j T}{n} = j h$, $0 \leq j \leq n $, this yields $R_{Y|V} = (\alpha_{ij})_{1 \leq i \leq d, 0 \leq j \leq n},$ with
\begin{itemize}
\item $\alpha_{ij} = \lambda_i^W \frac{2 e_i^W(t_j) - e_i^W(t_{j-1})- e_i^W(t_{j+1})}{h}$ for $j \notin \{0, n\}$,
\item $\alpha_{i0}= \lambda_i^W \Big( \left(e_i^W\right)'(t_0) - \frac{e_i^W(t_1)- e_i^W(t_0)}{h} \Big),$
\item $\alpha_{in} = \lambda_i^W \Big(\frac{e_i^W(t_n)- e_i^W(t_{n-1})}{h} - \left({e_i^W}\right)'(t_n) \Big).$
\end{itemize}
\vspace{4mm}
\par \noindent We now have a very fast and easy way to simulate the conditional distribution \eqref{eq:conditional_distribution} at our disposal.
\vspace{2mm}
\par
\begin{figure}[!ht]
	\begin{center}
	\psfrag{$0$}{$0$}
	\psfrag{$0.2$}{$0.2$}
	\psfrag{$0.4$}{$0.4$}
	\psfrag{$0.6$}{$0.6$}
	\psfrag{$0.8$}{$0.8$}
	\psfrag{$1$}{$1$}
	\psfrag{$-2$}{$-2$}
	\psfrag{$-1$}{$-1$}
	\psfrag{$2$}{$2$}
	\psfrag{$3$}{$3$}
	\psfrag{$-3$}{$-3$}
	\includegraphics[height=7cm]{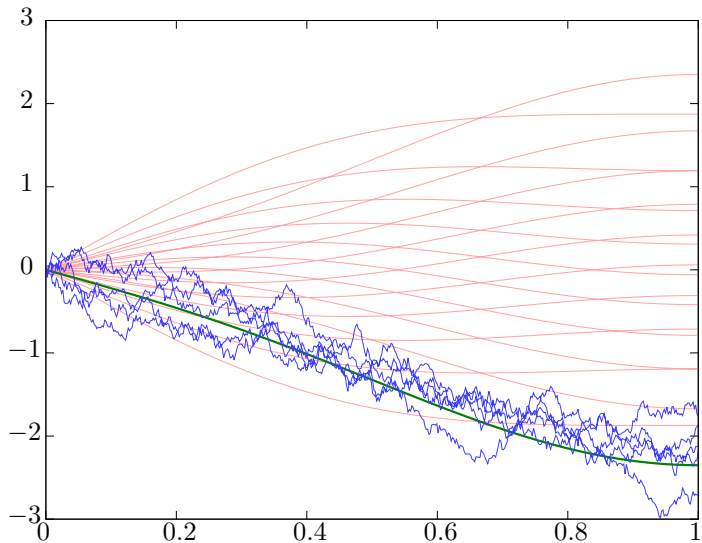}
	\caption{A few paths of the conditional distribution of Brownian motion,  given that its path belong to the $L^2$ Voronoi cell of the highlighted curve in the quantizer.}
	\label{fig:brownian_motion_conditional_simulation}
	\end{center}
\end{figure}

\begin{figure}[!ht]
\begin{minipage}[c]{.46\linewidth}
	\psfrag{0}{$0$}
	\psfrag{0.2}{$0.2$}
	\psfrag{0.4}{$0.4$}
	\psfrag{0.6}{$0.6$}
	\psfrag{0.8}{$0.8$}
	\psfrag{1}{$1$}
	\psfrag{-2.5}{$-2.5$}
	\psfrag{-1.5}{$-1.5$}
	\psfrag{-1}{$-1$}
	\psfrag{-0.5}{$-0.5$}
	\psfrag{0.5}{$0.5$}
	\psfrag{1.5}{$1.5$}
\includegraphics[height=4.5cm]{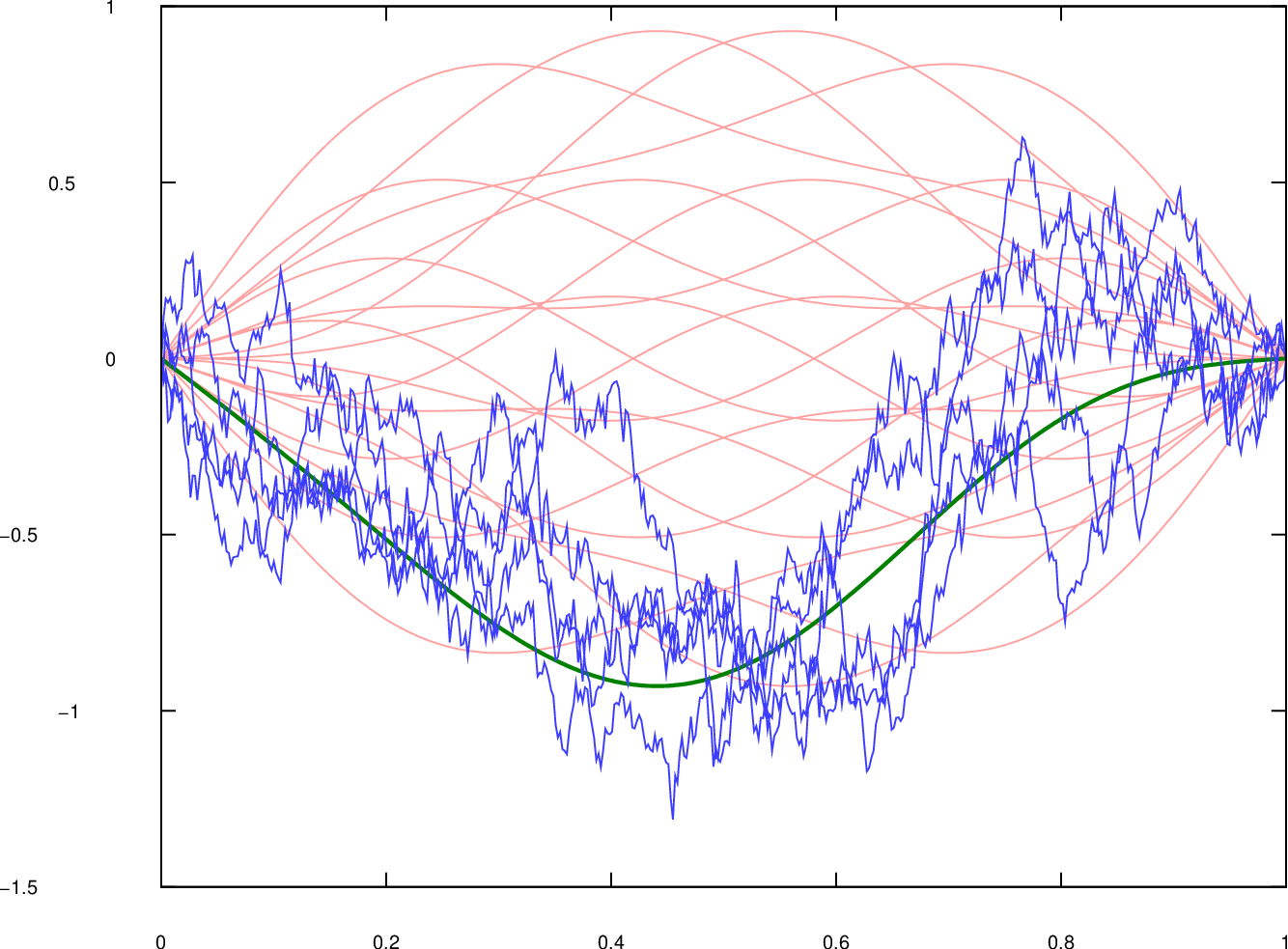}
\end{minipage} \hfill
\begin{minipage}[c]{.46\linewidth}
	\psfrag{-2.5}{$-2.5$}
	\psfrag{-2}{$-2$}
	\psfrag{-1.5}{$-1.5$}
	\psfrag{-1}{$-1$}
	\psfrag{-0.5}{$-0.5$}
	\psfrag{0}{$0$}
	\psfrag{0.5}{$0.5$}	
	\psfrag{1}{$1$}
	\psfrag{1.5}{$1.5$}
	\psfrag{2}{$2$}
	\psfrag{2.5}{$2.5$}	
	\psfrag{3}{$3$}
\includegraphics[height=4.5cm]{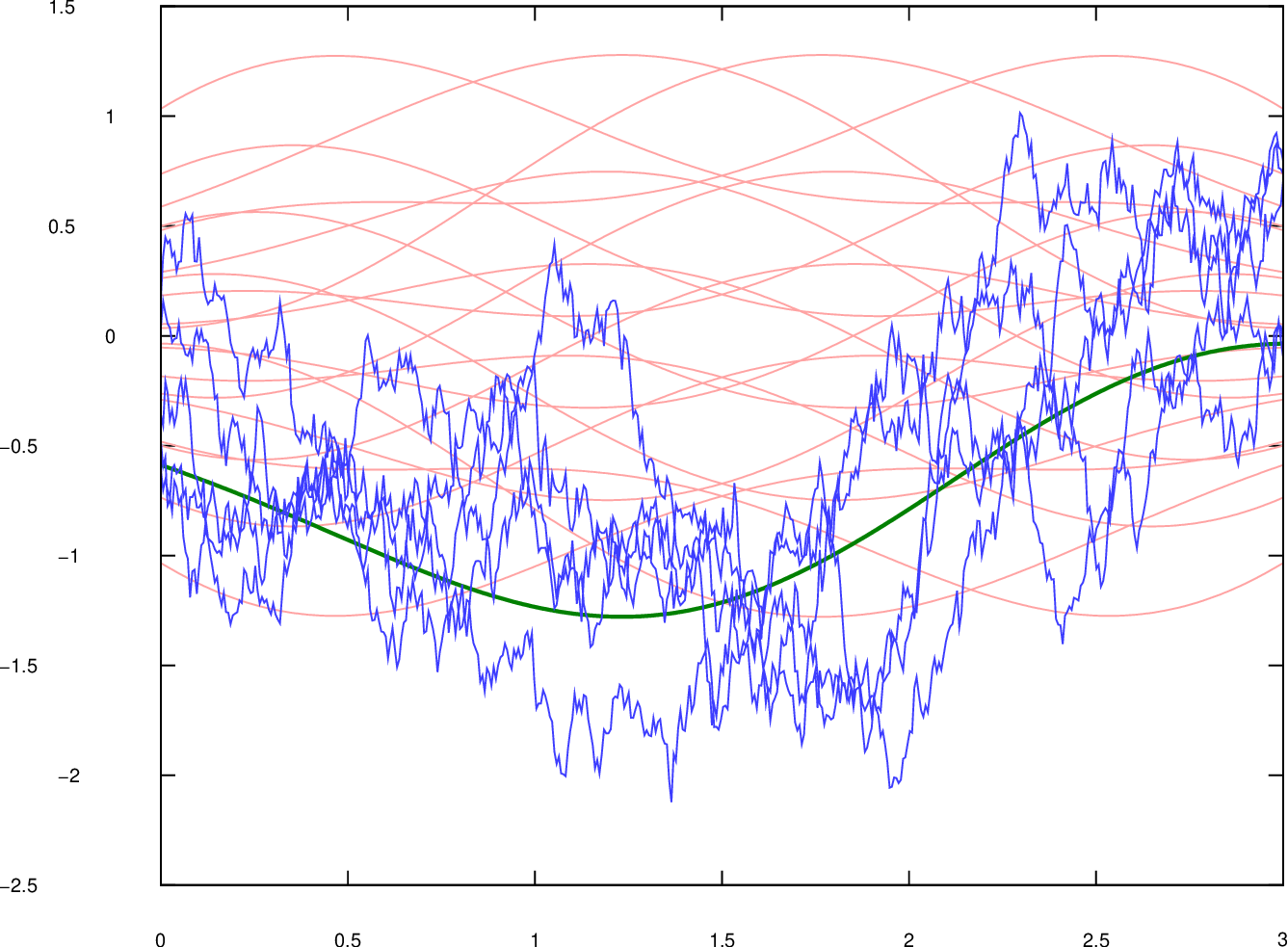}
\end{minipage}
\caption{A few paths of the conditional distribution of Brownian bridge (left) and the stationary Ornstein-Uhlenbeck process (right), given that its path belong to the $L^2$ Voronoi cell of the highlighted curve in the quantizer.}
\label{fig:other_conditional_simulation}
\end{figure}

 In Figures \ref{fig:brownian_motion_conditional_simulation} and \ref{fig:other_conditional_simulation}, we plot a few paths of the conditional distribution of various Gaussian processes given that they belong to a given $L^2$ Voronoi cell. The appearance of the drawing suggests to consider the method as a ``guided Monte Carlo simulation''. 
\subsection{Blind optimization procedures for the universal strata design}
\par We have seen in Section \ref{sec:quantization_based_stratif} that the quantity $d(\chi) = \Big(\sum\limits_{\chi_{\ui} \in \Gamma} p_{\ui} \sigma_{\ui} \Big) ^2$ is an upper bound of the variance of the estimator, given in Equation \eqref{eq:unbiased_stratif_estimator} in the case where the functional is $1$-Lipschitz continuous. 
Hence one may want to minimize this criterion instead of the $L^2$ quantization error. This yields the minimization problem
\begin{equation}\label{eq:schwartz_criteria1}
\mathcal{D}^{pq}_N := \min \Big\{ d(\chi ), \ \chi \in \mathcal{O}_{pq} (X,N) \Big\}
\end{equation}
instead of the minimization problem \eqref{eq:blind_optimization_dist}.
\par The same kind of blind optimization procedure as in Section \ref{sec:blind_optimization_quadratic} can be performed. Some values of the optimal decomposition for Brownian motion are given in Table \ref{tab:brownian_motion_optimal_decomposition_schwartz}. 
\begin{table}[!ht]
\begin{center}
\begin{tabular}{|c|c|c|c|}
\hline
$N$ & $N_{rec}$ & $d(\chi)$ & Product decomposition\\ 
\hline \hline
$1$ & $1$ & $0.5$ & $1$\\
$10$ & $10$	& $9.75689\cdot 10^{-2}$ & $5 \ \times \ 2$\\
$100$ & $96$ & $5.10548\cdot 10^{-2}$ & $12 \ \times \ 4 \ \times \ 2$\\
$1000$ & $966$ & $3.51289\cdot 10^{-2}$ & $23 \ \times \ 7 \ \times \ 3 \ \times \ 2$\\
$10000$ & $9984$ & $2.63721\cdot 10^{-2}$ & $26 \ \times \ 8 \ \times \ 4 \ \times \ 3 \ \times \ 2 \ \times \ 2$\\
\hline
\end{tabular}
\caption{Record of optimal product decomposition of Brownian motion with respect to the criterion \eqref{eq:schwartz_criteria1}. }
\label{tab:brownian_motion_optimal_decomposition_schwartz}
\end{center}
\end{table}
\par Optimal product decompositions for both Brownian bridge and Brownian motion and for a wide range of values of $N$ are available on the web site {\verb www.quantize.maths-fi.com } \cite{WebSiteGaussian} for download. When comparing decompositions for levels lesser $11000$, we notice that in the case of Brownian motion, the optimal decompositions for both criteria are ``almost'' always the same. The only values where decompositions differ are the ranges $270-271$ and $3328-3359$, and even then, the two criteria result in similar decompositions. Hence in practice, we can use the same database for both criteria. Nonetheless, in the case of Brownian bridge and Ornstein-Uhlenbeck processes, the optimal decompositions resulting from the two criteria differ more often. 
\subsection{Functional stratification of solutions of stochastic differential equations}
\par \noindent We consider the SDE 
\begin{equation}\label{eq:oned_diffusion}
dF_t = b(t,F_t) dt + \sigma(t,F_t) dX_t, \hspace{3mm} t \in [0,T], \ \quad F_0 = f_0
\end{equation}
where $X$ is a centered continuous Gaussian semimartingale starting from $0$ and where $b$ and $\sigma$ are Borel functions, Lipschitz continuous in $x$ uniformly in $t \in [0,T]$ such that $|b(\cdot,0)|+|\sigma(\cdot,0)|$ is bounded over $[0,T]$. In this situation, \eqref{eq:oned_diffusion} admits a unique strong solution $X$ and $\sup\limits_{t \in [0,T]} |X_t|$ has $r$-moments for every $r \in (0,\infty)$. 
\begin{remark}
\par In this case, thanks to Fernique's theorem, the continuity assumption on the Gaussian process ensures that $\int_0^T \E[X_s^2] ds < \infty$ and the continuity of the covariance function, (see \cite[VIII.3]{JansonGaussianHilbertSpaces}). 
\end{remark}
\par The most common approach to perform a Monte Carlo simulation with the solution of such a stochastic differential equation, is to use a discretization scheme like the Euler scheme \cite{GlassermanMonteCarlo}. In this setting, we propose to simply replace Gaussian process $X$ by a stratified version of $X$ in the Euler scheme. This approach is justified in many ways: 
\begin{enumerate}
\item In \cite{CorlayPartialQuantization}, using filtration enlargement techniques, it is proved that under some additional hypothesis on the Gaussian semimartingale $X$, its conditional distribution in a strata is still a semimartingale with respect to its own filtration. This additional hypothesis is satisfied by Brownian motion, Brownian bridge and Ornstein-Uhlenbeck processes. Therefore, plugging the stratified Euler scheme into the SDE amounts to using the Euler scheme of these conditional stochastic differential equations. 
\item In the one-dimensional setting, if we make the additional hypothesis that $\sigma \in C^1([0,T]\times \R,\R)$ is positive and bounded, as soon as the drift of the Lamperti transform of the SDE \eqref{eq:oned_diffusion} is Lipschitz continuous, the unique strong solution of \eqref{eq:oned_diffusion}, seen as a functional of the underlying Gaussian process $X$ is $\|\cdot\|_p$-Lipschitz continuous \cite{LuschgyPagesFunctional5}. Hence we can apply the results of Section \ref{sec:quantization_based_stratif} on universal stratification for Lipschitz continuous functionals. 
\item The function $\left(X_{t_0}, \ldots, X_{t_n}\right) \mapsto \left(X_{t_1}-X_{t_0}, \ldots, X_{t_n}-X_{t_{n-1}}\right)$ that maps the marginals of Brownian motion to the corresponding increments used in the Euler scheme, is a linear map from $\R^{n+1}$ to $\R^n$ and thus Lipschitz continuous as well. 
\end{enumerate}
\section{Application to option pricing}\label{sec:option_pricing}
\par The special case of Brownian motion allows us to use functional stratification as a generic variance reduction method for the case of functionals of Brownian diffusions, even in the multidimensional case, regardless of how the Brownian paths are correlated or used afterwards, to drive the diffusion of an underlying stock, a stochastic volatility process or a discount factor. As it only impacts how the Monte Carlo simulation is input with Brownian paths,  our approach is easier to implement in a practical setting than adaptive variance reduction methods, which generally require a control loop. 
\par In this section, we study the performance of our method in simple one-dimensional cases. We begin with the case of a continuous-time Up-In Call option in the Black and Scholes model, for which a closed-form expression is known, and used as a Benchmark. 
\subsection{Benchmark with an Up-In Call in the Black and Scholes model}\label{sec:UIC_Benchmark}
\par We evaluate our method in the case of a path dependent option where a reference value can easily be computed: an Up-In Call barrier option in the Black and Scholes model. For the sake of simplicity, we assume that there is no drift (no interest rate and no dividend). There is a closed-form expression for the continuous barrier option price, but we must resort to a numerical approximation \cite{GlassContinuousCorrection} (yet very accurate) on that closed-form expression to get the price in the case of discrete dates for the barrier. The total size of the Monte Carlo sample is $100 000$ in every case. 
\par We price the Up-In Call option with different values of the initial spot $S$, the strike $K$, the barrier $H$, the volatility $\sigma$, the maturity $T$, and the number of fixing dates for the discrete barrier $n$. In every case, a $95 \%$ confidence interval is given. So is the variance of the estimator. 
\par The numerical results are reported in Table \ref{fig:var_reduc_20_points} when using the method with $20$ stratas and Table \ref{fig:var_reduc_100_points} when using the method with $100$ stratas. In this tables, the first column correspond to Broadie and Glasserman's closed-form expression proxy. The second one corresponds to a simple Monte Carlo estimator. The last three columns correspond to stratified sampling estimators with different simulation allocation strategies. 
\par The ``natural weights'' column stands for the allocation budget of Equation \eqref{eq:sub_optimal_choice}. The ``Lip.-optimal weights'' column stand for the ``universal stratification'' budget allocation proposed in Section \ref{sec:quantization_based_stratif}. In these two cases, we have an explicit allocation rule which does not depend on the payoff function. The last column, ``optimal weights'' corresponds to an estimation of the optimal budget allocation given in Equation \eqref{eq:optimal_budget}. 

\begin{table}[!ht]
	\begin{center}
	{
	\scriptsize
	\begin{tabular}{|c|c|c|c|c|c|}
	\hline
	Parameter				& Broadie \&		&	Simple 		&	Strat. estimator	& Strat. estimator		& Strat. estimator	\tabularnewline
        Values		& Glasserman's		&	estimator	&	natural weights	& Lip.-optimal weights	& optimal weights	\tabularnewline	
							& proxy				&				&						&						&					\tabularnewline
	\hline
	\hline
	$S = 100$, $K = 100$	&			& $14.0379$				& $13.9281$				& $13.9283$				& $13.9364$				\tabularnewline	
	$H=125$, $\sigma=0.3$,	& $13.9597$	& $[13.8705,14.2053]$	& $[13.8491,14.0071]$	& $[13.8519,14.0047]$	& $[13.8827,13.9901]$	\tabularnewline
	$T=1.5$, $n=365$ 		& 			& $\var = 729.2518$		& $\var = 162.4650$		& $\var = 151.9481$		& $\var = 75.1319$		\tabularnewline
	\hline
	$S = 100$, $K = 100$	&			& $1.4206$			& $1.3659$			& $1.3510$			& $1.3602$			\tabularnewline	
	$H=200$, $\sigma=0.3$,	& $1.3665$	& $[1.3442,1.4969]$	& $[1.3106,1.4211]$	& $[1.3039,1.3981]$	& $[1.3472,1.3732]$	\tabularnewline
	$T=1$, $n=365$ 			& 			& $\var = 151.6366$	& $\var = 79.5118$	& $\var = 57.7425$	& $\var = 4.4053$	\tabularnewline
	\hline
	\end{tabular}
	}
	\caption{Numerical results for the Up-In Call option, with $20$ stratas. }
	\label{fig:var_reduc_20_points}
	\end{center}
\end{table}
\begin{table}[!ht]
	\begin{center}
	{
	\scriptsize
	\begin{tabular}{|c|c|c|c|c|c|}
	\hline
	Parameter				& Broadie \&	&	Simple 		&	Strat. estimator	& Strat. estimator		& Strat. estimator	\tabularnewline
	Values						& Glasserman's	&	estimator	&	natural weights	& Lip.-optimal weights	& optimal weights	\tabularnewline
							& proxy			&				&						& 						& 					\tabularnewline
	\hline
	\hline
	$S = 100$, $K = 100$	&			& $14.0379$				& $13.9382$				& $13.9511$				& $13.9483$				\tabularnewline	
	$H=125$, $\sigma=0.3$,	& $13.9597$	& $[13.8705,14.2053]$	& $[13.8720,14.0043]$	& $[13.8874,14.0150]$	& $[13.9047,13.9919]$	\tabularnewline
	$T=1.5$, $n=365$ 		& 			& $\var = 729.2518$		& $\var = 114.0634$		& $\var = 105.8760$		& $\var = 49.5071$		\tabularnewline
	\hline
	$S = 100$, $K = 100$	& 			& $1.4206$			& $1.3296$			& $1.3493$			& $1.3611$ 					\tabularnewline	
	$H=200$, $\sigma=0.3$,	& $1.3665$	& $[1.3442,1.4969]$	& $[1.2825,1.3768]$	& $[1.3093,1.3893]$	& $[1.3508,1.3715]$			\tabularnewline
	$T=1$, $n=365$ 			& 			& $\var = 151.6366$	& $\var = 57.8899$	& $\var = 41.6666$ & $\var = 2.8099$ 			\tabularnewline
	\hline
	\end{tabular}
	}
	\caption{Numerical results for the Up-In Call option, with $100$ stratas. }
	\label{fig:var_reduc_100_points}
	\end{center}
\end{table}
\subsection{Test with an Auto-Call in the CEV model}\label{sec:autocall_CEV}
\par We assume that the stock follows a CEV model with no drift $dS_t = \sigma S_t^{\frac{\beta}{2}} dW_t, \hspace{4mm} 0 < \beta < 2$. We used the Euler scheme on $\ln(S_t)$, which satisfies the SDE $d\ln(S_t) = -\frac{\sigma^2}{2} S_t^{\beta-2} dt + \sigma S_t^{\frac{\beta}{2}-1} dW_t$.
\vspace{2mm}
\par \noindent \textbf{Description of the Auto-Call payoff:} 
\par Let $S_t$ be the stock price and $0 = t_0 < t_1 < \cdots < t_n = T$ be the observation dates. $K$ and $H$ are the ``strike'' and the ``barrier'' values. $P$ denotes the ``nominal'', and $C$ a zero-coupon bond of maturity $T$. 
\par At the first date $t_1$ of the schedule, if $S_{t_1}>K$, the holder of the option receives $(1+C) P$ and the contract expires. If $S_{t_1}\leq K$, he waits until the second date of the schedule. If $S_{t_2}>K$, the holder gets $(1+C) P$ and the contract expires. And so on... If $S_t$ does not reach $K$ on $[0,T)$, the contract is exerciced as follows: if $S_T>K$, the holder gets $(1+C) P$. If $H < S_T \leq K$, the holder gets $P$ and if $S_T \leq H$, he gets $P \frac{S_T}{K}$.
\par The numerical results are reported in Table \ref{fig:autocall_var_reduc_20_points} when implementing the stratification method with $20$ and $50$ stratas. The parameters of the model are $\beta = 1.5$, $S_0 = 100$, $\sigma = 0.3$. For the payoff, $K = 110$, $H = 80$, $P=100$, $C=0.07$. The considered observation dates are $\{1,2,3\}$. The number of time steps in the Euler scheme is $300$ and the total size of the Monte Carlo sample is $100 000$ in every case. 
\begin{table}[!ht]
	\begin{center}
	{
	\scriptsize
	\begin{tabular}{|c|c|c|c|c|}
	\hline
	Number of strata	& Simple 				&	Strat. estimator	& Strat. estimator		& Strat. estimator		\tabularnewline
						& estimator				&	natural weights	& Lip.-optimal weights	& optimal weights		\tabularnewline
	\hline 
	\hline
						& $99.0598$				& $99.0839$				& $99.0886$				& $99.0477$				\tabularnewline	
	$20$				& $[98.9887,99.1310]$	& $[99.0438,99.1239]$	& $[99.0488,99.1284]$	& $[99.0184,99.0769]$	\tabularnewline
						& $\var = 131.8089$		& $\var = 41.8067$		& $\var = 41.2888$		& $\var = 22.2549$		\tabularnewline
	\hline
						& $99.0598$				& $99.0507$				& $99.0790$				& $99.0444$				\tabularnewline	
	$50$				& $[98.9887,99.1310]$	& $[99.0129,99.0886]$	& $[99.0414,99.1166]$	& $[99.0179,99.0709]$	\tabularnewline
						& $\var = 131.8089$		& $\var = 37.3150$		& $\var = 36.8408$		& $\var = 18.2954$		\tabularnewline
	\hline
	\end{tabular}
	}
	\caption{Numerical results for the Auto-Call option in the CEV model, with $20$ and $50$ stratas. }
	\label{fig:autocall_var_reduc_20_points}
	\end{center}
\end{table}
\subsection{Test with an Asian straddle in the one-factor Schwartz model}\label{sec:asian_schwartz}
\par Here, we stand in the case of a stock which follows the following SDE:
\begin{equation}
dS_t = \theta (\alpha - \ln S_t) S_t dt + \sigma S_t dW_t,
\end{equation}
under the risk-neutral probability. The stochastic process $X = \ln(S)$ is an Ornstein-Uhlenbeck process:
\begin{equation}
dX_t = \theta (\mu- X_t) dt + \sigma dW_t \hspace{5mm} \textnormal{with } \mu = \alpha - \frac{\sigma^2}{2 \theta}.
\end{equation}
\par This model was proposed by Schwartz in \cite{SchwartzCommo97}. Such exponentials of Ornstein-Uhlenbeck processes are commonly met in commodity derivatives. One particularity in these markets is that the spot is generally not directly traded. Therefore, the underlyings of derivatives are generaly futures. Still, we use this one-factor ``toy'' model as a simple case study for our variance reduction method. 
\par The considered payoff is an \textit{Asian straddle} option on a discrete schedule of observation dates $t_0 < \cdots < t_n = T$. $K$ is the ``strike'' of the option whose payoff is $\left| \frac{1}{n+1} \sum\limits_{k=0}^n S_{t_k} -K \right|$.
\par We perform a functional stratified sampling of the Ornstein-Uhlenbeck process. Optimal product decompositions for the criterion \eqref{eq:schwartz_criteria1} are used and available in Table \ref{fig:schwartz_var_reduc_20_points} where the numerical results are reported.
\par The parameters are $S_0 = 100$, $\theta = 0.3$, $\alpha = \ln(110)$, $\sigma = 0.3$ and $K = 100$. The total size of the Monte Carlo sample is $100 000$ in every case. The observation dates are $\left(i \frac{T}{n}\right)_{i=\{0, \ldots, n\}}$ with $T=3$ and $n=36$. 
\begin{table}[!ht]
	\begin{center}
	{
	\scriptsize
	\begin{tabular}{|c|c|c|c|c|}
	\hline
	Number of strata			& Simple 			&	Strat. estimator	& Strat. estimator		& Strat. estimator		\tabularnewline
	and product decomposition	& estimator		&	natural weights	& Lip.-optimal weights	& optimal weights		\tabularnewline
	\hline 
	\hline
								& $17.5393$				& $17.6140$				& $17.6118$				& $17.6240$			\tabularnewline	
	$20$						& $[17.4504,17.6282]$	& $[17.5871,17.6408]$	& $[17.5856,17.6378]$	& $[17.6006,17.6477]$\tabularnewline
	$20 = 10 \times 2$			& $\var = 205.9375$		& $\var = 18.8041$		& $\var = 17.5502$		& $\var = 14.6363$	\tabularnewline
	\hline
								& $17.5393$				& $17.6101$				& $17.6122$				& $17.6147$			\tabularnewline	
	$100$						& $[17.4504,17.6282]$	& $[17.5850,17.6351]$	& $[17.5884,17.6360]$	& $[17.5932,17.6362]$\tabularnewline
	$100 = 10 \times 5 \times 2$& $\var = 205.9375$		& $\var = 16.2945$		& $\var = 14.7316$		& $\var = 12.0112$	\tabularnewline
	\hline
	\end{tabular}
	}
	\caption{Numerical results for the Asian straddle option in Schwartz's model, with $20$ and $100$ stratas. }
	\label{fig:schwartz_var_reduc_20_points}
	\end{center}
\end{table}
\par To perform this computation, one needs to use a non-centered Ornstein-Uhlenbeck quantizer. Building such a quantizer is a straightforward extension of the centered case. As showed in Section \ref{sec:karhunen_loeve_basis}, if $X$ is an Ornstein-Uhlenbeck process on $[0,T]$ following the dynamic
$dX_t = \theta (\mu-X_t) dt + \sigma dW_t$, $X_0 \simdist \mathcal{N}(m_0,\sigma_0^2)$, with nonzero values of $\mu$ and $m_0$, we have 
\begin{equation}
X_t = \underbrace{m_0 e^{-\theta t} + \mu (1-e^{- \theta t})}_{(1) =\textnormal{non-stochastic path}} + \left(\begin{array}{ccc}\textnormal{centered Ornstein-Uhlenbeck process}\\ \textnormal{corresponding to } m_0 = \mu = 0 \end{array}\right). 
\end{equation}
\par
\begin{figure}[!ht]
	\begin{center}
	\psfrag{$0$}{$0$}
	\psfrag{$0.5$}{$0.5$}
	\psfrag{$1$}{$1$}
	\psfrag{$1.5$}{$1.5$}
	\psfrag{$2$}{$2$}
	\psfrag{$2.5$}{$2.5$}
	\psfrag{$3$}{$3$}
	\psfrag{$4.6$}{$4.6$}
	\psfrag{$4.8$}{$4.8$}
	\psfrag{$5$}{$5$}
	\psfrag{$5.2$}{$5.2$}
	\psfrag{$5.4$}{$5.4$}
	\psfrag{$5.6$}{$5.6$}
	\psfrag{$5.8$}{$5.8$}
	\psfrag{$6$}{$6$}
	\includegraphics[width=8cm]{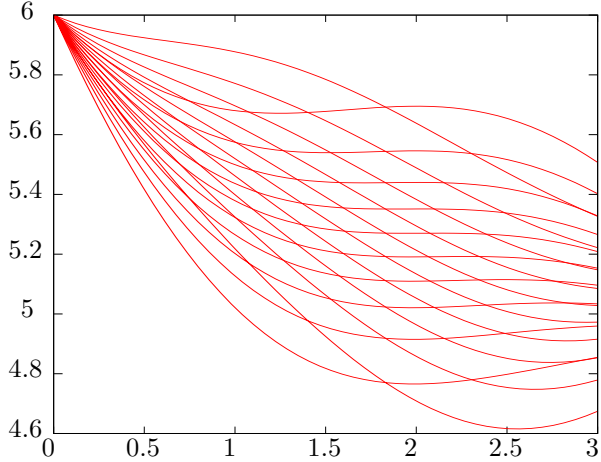}
	\caption{Functional $10 \times 2$-product quantizer of an Ornstein-Uhlenbeck process starting from $X_0 = 6$ defined by the diffusion $dX_t = \theta(\mu - X_t)dt + \sigma dW_t$ with $\mu = 5$, $\sigma = 0.3$ and $\theta = 0.8$ on $[0,3]$.} 
	\label{fig:non_central_ou}
	\end{center}
\end{figure}

 \noindent In Figure \ref{fig:non_central_ou}, we display the functional product quantizer of a non-centered Ornstein-Uhlenbeck process. 

\subsection{Comments on the numerical results}
\par ``Lipschitz-optimal'' strata and weights are not more difficult to compute than the ``natural'' scheme \eqref{eq:sub_optimal_choice} since all the involved parameters are known. This stratified sampling method does not depend on the payoff function but only on the distribution of the underlying asset which means that it can be plugged upstream in the Monte Carlo simulator. In terms of variance reduction, universal stratification is all the more preforming as the variance of the plain estimator is high, most likely because more strata are ``activated''. When the payoff function is symmetric, like with the \emph{Asian straddle} it achieves up to $90 \%$ of the variance reduction achieved by a payoff-dependent dedicated method like the one devised in \cite{JourdainStratification1}.
\begin{appendices}
\section{The Karhunen-Loève expansion of Ornstein-Uhlenbeck processes}\label{sec:karhunen_loeve_basis}
\par In this section, we derive the Karhunen-Loève expansion of the Ornstein-Uhlenbeck process. Proposition \ref{prop:OU_KL} brings the results together. Section \ref{sec:OU_numerical} presents the numerical method for computing this expansion. 
\subsection{The Ornstein-Uhlenbeck process}
The Ornstein-Uhlenbeck process is defined by the SDE
\begin{equation}\label{eq:ornstein_uhlenbeck_sde}
d X_t = \theta (\mu - X_t) dt + \sigma dW_t, \hspace{8mm} \textnormal{ with } \sigma \geq 0 \textnormal{ and } \theta > 0,
\end{equation}
which gives
\begin{equation}{\label{eq:ornstein_uhlenbeck_solution}}
X_t = X_0 e^{-\theta t} + \mu (1-e^{-\theta t}) + \int_0^t \sigma e^{\theta (s-t)} dW_s.
\end{equation}
\par \noindent We assume that $X_0$ is Gaussian ($X_0 \simdist \mathcal{N}(m_0,\sigma^2_0)$) and independent from $W$. We have $\E[X_t] = m_0 e^{-\theta t} + \mu (1-e^{-\theta t})$
and $\cov(X_s,X_t) = \frac{\sigma^2}{2\theta} e^{-\theta (s+t)} \left(e^{2 \theta \min(s,t)} - 1 \right) + \sigma_0^2 e^{-\theta (s+t)}.$ Moreover $\lim\limits_{t \to \infty} \var(X_t) = \frac{\sigma^2}{2 \theta}$ (the long-term variance). If the initial variance $\sigma_0^2$ is equal to long-term variance $\frac{\sigma^2}{2 \theta}$, $X$ is stationary and the covariance function is given by $\cov(X_s,X_t) = \frac{\sigma^2}{2\theta} e^{-\theta |s-t|}$. The total variance of the process on $[0,T]$ is
$$
\|X\|_2^2 = \int_0^T \var(X_s) ds = \frac{\sigma^2 T}{2\theta} + \left(\sigma_0^2 - \frac{\sigma^2}{2 \theta}\right) \left(\frac{1}{2 \theta} - \frac{e^{-2 \theta T}}{2 \theta} \right).
$$
\subsection{The Ornstein-Uhlenbeck covariance operator}
\par The Ornstein-Uhlenbeck covariance operator is given by
\begin{equation}\label{ornstein_uhlenbeck_covariance_operator}
T^{OU} f(t) = \int_0^T \frac{\sigma^2}{2 \theta} e^{-\theta (s+t)} \left( e^{2 \theta \min(s,t)} -1 \right) f(s) ds + \int_0^T \sigma_0^2 e^{-\theta (s+t) } f(s) ds. 
\end{equation}
\par \noindent \textbf{Computing the Karhunen-Loève expansion of the Ornstein-Uhlenbeck process}
\par \noindent $T^{OU}$ is a compact Hermitian positive operator on the separable Hilbert space $L^2([0,T])$. Hence there exists an orthonormal basis consisting of eigenvectors of $T^{OU}$ and eigenvalues are real and nonnegative. Moreover $\left\|T^{OU}\right\|^2 \leq \frac{\sigma^2 T}{2 \theta} + \frac{\sigma^2}{4 \theta^2} \left( e^{-2 \theta T} - 1 \right) $. We have 
{\small
$$
T^{OU} f(t) = \int_0^t \frac{\sigma^2}{2 \theta} e^{\theta (s-t)} f(s) ds + \int_t^T \frac{\sigma^2}{2 \theta} e^{\theta (t-s)} f(s) ds + \int_0^T \left(\sigma_0^2 - \frac{\sigma^2}{2 \theta} \right) e^{- \theta(s+t)} f(s) ds.
$$
}
\begin{prop}
\par If $f\in C([0,1])$, and if $g = T^{OU}f$, then
	\begin{equation}\label{eq:ornstein_uhlenbeck_intermediate1}
	g''-\theta^2 g = -\sigma^2 f,
	\end{equation}
	with
	\begin{equation}\label{eq:ornstein_uhlenbeck_intermediate2}
	\sigma_0^2 g'(0) = \left( \sigma^2 - \theta \sigma_0^2 \right) g(0) \hspace{8mm} \textnormal{and} \hspace{8mm} g'(T) = -\theta g(T).
	\end{equation}
\end{prop}
\noindent \textbf{Proof:}
{\small
$$
g(t) = \int_0^t \frac{\sigma^2}{2 \theta} e^{\theta (s-t)} f(s) ds + \int_t^T \frac{\sigma^2}{2 \theta} e^{\theta (t-s)} f(s) ds + \int_0^T \left(\sigma_0^2 - \frac{\sigma^2}{2 \theta}\right) e^{- \theta(s+t)} f(s) ds.
$$
$$
g'(t) = -\frac{\sigma^2}{2} \int_0^t e^{\theta (s-t)} f(s) ds + \frac{\sigma^2}{2} \int_t^T e^{\theta (t-s)} f(s) ds -\left(\theta \sigma_0^2 - \frac{\sigma^2}{2 }\right) \int_0^T e^{-\theta (s+t)} f(s) ds \\
$$
$$
g''(t) = \frac{\sigma^2 \theta}{2} \left( \int_0^t f(s) e^{\theta(s-t)} ds+ \int_t^T f(s) e^{\theta (t-s)}ds \right) + \theta \int_0^T \left(\theta \sigma_0^2 - \frac{\sigma^2}{2}\right) e^{-\theta (s+t)} f(s) ds - \sigma^2 f(t).
$$
\par \noindent we get $g''(t) = \theta^2 g(t) - \sigma^2 f(t).$ Moreover, Equation \eqref{eq:ornstein_uhlenbeck_intermediate2} comes when identifying expressions with $t=0$ and $t=T$.
}
\myqed
\begin{prop}
\par Conversely, if $g\in C^2([0,T])$ and if functions $f$ and $g$ satisfy Equations \eqref{eq:ornstein_uhlenbeck_intermediate1} and \eqref{eq:ornstein_uhlenbeck_intermediate2} then $g = T^{OU}f$. 
\end{prop}
\noindent \textbf{Proof:} Computing $T^{OU}g''$ yields: 
$$
T^{OU}g''(t) = \int_0^t \frac{\sigma^2}{2 \theta} e^{\theta (s-t)} g''(s) ds + \int_t^T \frac{\sigma^2}{2\theta} e^{\theta (t-s)} g''(s) ds + \int_0^T \left(\sigma_0^2 - \frac{\sigma^2}{2\theta}\right) e^{-\theta (s+t)}g''(s) ds.
$$
An integration by parts yields
$$
\begin{array}{lll}
T^{OU}g'' 	&= -\sigma_0^2 g'(0) e^{-\theta t} - \sigma^2 g(t) + \frac{\sigma^2}{2} g(0) e^{- \theta t} - \left(\theta \sigma_0^2 -\frac{\sigma^2}{2} \right) g(0) e^{-\theta t} + \theta^2 T^{OU}g(t) \\
			&= -\sigma^2 g(t) + \theta^2 T^{OU}g(t) \hspace{4mm} \textnormal{thanks to \eqref{eq:ornstein_uhlenbeck_intermediate2}.}
\end{array}
$$
\myqed
\par \noindent Now, by necessary conditions, $T^{OU}f = \lambda f \Leftrightarrow \sigma^2 g = \lambda (\theta^2 g - g'')$. We obtain
\begin{equation}\label{eq:ornstein_uhlenbeck_ode}
\lambda g'' + (\sigma^2-\lambda \theta^2)g =0.
\end{equation}
\par \noindent Hence the solution of the ordinary differential equation \eqref{eq:ornstein_uhlenbeck_ode} on $[0,T]$ has the form
$g(t) = A \cos(\omega t) + B \sin(\omega t),$
with $\omega = \sqrt{\frac{\sigma^2-\lambda \theta^2}{\lambda}} \Leftrightarrow \lambda = \frac{\sigma^2}{\omega ^2 + \theta^2}$. Equation \eqref{eq:ornstein_uhlenbeck_intermediate2} yields $\omega B \sigma_0^2 = (\sigma^2-\theta \sigma_0^2) A$. Hence, we have $g(t) = K \left( \omega \sigma_0^2 \cos(\omega t) + (\sigma^2-\theta \sigma_0^2) \sin(\omega t) \right)$, so that Equality \eqref{eq:ornstein_uhlenbeck_intermediate2} yields
\begin{equation}\label{eq:ornstein_uhlenbeck_functional_equation}
\omega \sigma^2 \cos(\omega T) + \left( -\omega^2 \sigma_0^2 + \theta \sigma^2 - \theta^2 \sigma_0^2 \right) \sin(\omega T) = 0.
\end{equation}
\par \noindent Conversely, the same calculation shows that $\lambda_n \in \left]0,\left\|T^{OU}\right\|_2\right]$ is an eigenvalue of $T^{OU}$ if and only if Equality \eqref{eq:ornstein_uhlenbeck_functional_equation} holds.
\begin{prop}\label{prop:OU_KL} Finally, if $(\omega_n)_{n \geq 1}$ is the increasingly sorted sequence of the positive solutions of \eqref{eq:ornstein_uhlenbeck_functional_equation}, the Karhunen-Loève eigensystem $\left(e_n^{OU},\lambda_n^{OU}\right)_{n \geq 1}$ of the Ornstein-Uhlenbeck process is
\begin{itemize}
\item $\lambda^{OU}_n = \frac{\sigma^2}{\omega_n^2 + \theta^2},$ and 
\item $e_n^{OU}(t) = K_n \left( \omega_n \sigma_0^2 \cos(\omega_n t) + (\sigma^2 - \theta \sigma_0^2) \sin(\omega_n t)\right)$
for $n \geq 1$, where $K_n$ is the normalization constant. 
\end{itemize}
\par \noindent If $(\sigma,\sigma_0) \neq (0,0)$, $K_n$ is given by
\begin{multline}
1/K_n^2 = \frac{1}{2 \omega_n} \sigma_0^2 (\sigma^2 - \theta \sigma_0^2) \left(1 - \cos(2 \omega_n T)\right) + \frac{1}{2}\sigma_0^4 \omega_n^2 \left(T + \frac{1}{2 \omega_n} \sin(2 \omega_n T) \right) \\
	+ \frac{1}{2}(\sigma^2-\theta \sigma_0^2)^2\Big(T-\frac{1}{2\omega_n}\sin(2 \omega_n T) \Big).
\end{multline}
\end{prop}
\par \noindent \textbf{Case of a deterministic starting point:} In this case ($\sigma_0 = 0$), we have 
$$
e_n^{OU}(t) = \frac{1}{\sqrt{\frac{T}{2}-\frac{\sin(2 \omega_n T)}{4 \omega_n}}} \sin(\omega_n t).
$$
\par \noindent \textbf{Stationary case:} In the stationary case $(\sigma_0^2 = \frac{\sigma^2}{2 \theta})$, we have $e_n^{OU}(t) = C_n \Big( \omega_n \cos(\omega_n t) + \theta \sin(\omega_n t)\Big)$, where $C_n$ is the normalization constant. $C_n$ is given by
$$
1/C_n^2 = \frac{\theta}{2} \Big(1-\cos(2\omega_n T)\Big)+ \frac{\omega_n^2}{2} \Big( T + \frac{\sin(2 \omega_n T)}{2 \omega_n}\Big) + \frac{\theta^2}{2}\Big( T-\frac{\sin(2 \omega_n T)}{2 \omega_n}\Big).
$$
\subsection{Numerical computation of the Karhunen-Loève expansion of the Ornstein-Uhlenbeck process}\label{sec:OU_numerical}
\par This section focuses on the computation of the positive solutions to \eqref{eq:ornstein_uhlenbeck_functional_equation}. 
\subsubsection{Deterministic starting point}
\par In this case ($\sigma_0 =0$), we can check that elements of $\left\{ \frac{\pi}{2T} + k \frac{\pi}{T}, \ k \in \N \right\}$ are not solutions of Equation \eqref{eq:ornstein_uhlenbeck_functional_equation}. As a consequence, the equation comes to
\begin{equation}\label{eq:conditional_ornstein_uhlenbeck_functional_equation}
\theta \tan(\omega T) = - \omega.
\end{equation}
\par
\begin{figure}[!ht]
	\begin{center}
	\psfrag{w1}{$\omega_1$}
	\psfrag{w2}{$\omega_2$}
	\psfrag{w3}{$\omega_3$}
	\psfrag{w4}{$\omega_4$}
	\psfrag{$0$}{$0$}
	\psfrag{tan}{ \small{$\tan(\omega T)$}}
	\psfrag{affine}{\small{$-\omega/\theta$}}
	\psfrag{pi/(2T)}{$\frac{\pi}{2T}$}
	\psfrag{pi/(2T)+pi/T}{$\frac{\pi}{2T} + \frac{\pi}{T}$}
	\psfrag{pi/(2T)+2pi/T}{$\frac{\pi}{2T} + \frac{2\pi}{T}$}
	\psfrag{pi/(2T)+3pi/T}{$\frac{\pi}{2T} + \frac{3\pi}{T}$}
	\includegraphics[width=8cm]{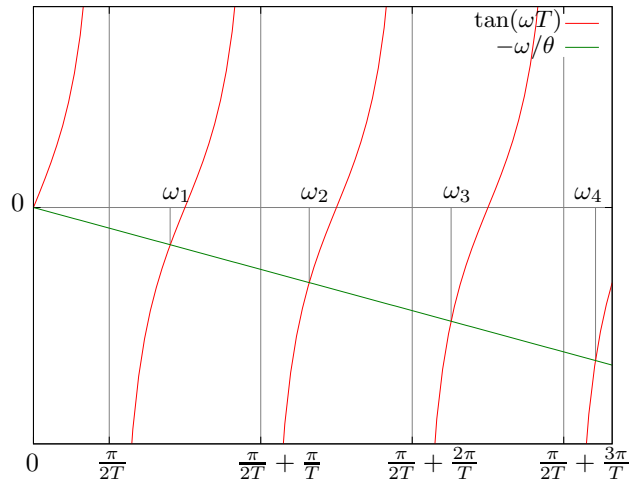}
	\caption{\textbf{(Deterministic starting point)}. Solutions of \eqref{eq:conditional_ornstein_uhlenbeck_functional_equation}. (Ornstein-Uhlenbeck process starting from a fixed point $X_0$, $\sigma_0 =0$.) Parameter values are $T=3$, $\sigma=1$ and $\theta=3$.}
	\label{fig:conditional_ornstein_uhlenbeck_eigen}
	\end{center}
\end{figure}

 \noindent The case where $\theta = 0$ comes to the case of Brownian motion, hence we assume that $\theta \neq 0$. 
Solutions of this equation are illustrated in Figure \ref{fig:conditional_ornstein_uhlenbeck_eigen}. 
We can show that there is a unique solution $\omega_n$ in each interval $\left(\frac{n\pi}{T}-\frac{\pi}{2T},\frac{n\pi}{T} \right)$,
for $n \in \{ 1, 2, \ldots \}$ and that $\lim\limits_{n\to \infty}\omega_n - \Big(\frac{n\pi}{T}-\frac{\pi}{2T}\Big) = 0$. 

\subsubsection{Non-deterministic starting point}
\par Let us assume now that $\sigma_0 \neq 0$ and consider Equation \eqref{eq:ornstein_uhlenbeck_functional_equation} again. The term $-\omega^2 \sigma_0^2 + \theta \sigma^2 - \theta^2 \sigma_0^2$ never vanishes on $\left(0,+\infty\right)$ if $\theta^2 \sigma_0^2-\theta \sigma^2 \geq 0$.
\vspace{1mm}
\par \noindent \textbf{First case: $\theta^2 \sigma_0^2-\theta \sigma^2 \geq 0$.} In this case \eqref{eq:ornstein_uhlenbeck_functional_equation} gives
\begin{equation}\label{eq:general_ornstein_uhlenbeck_functional_equation}
\tan(\omega T) = \frac{\omega \sigma^2}{\omega^2 \sigma_0^2+ \theta^2 \sigma_0^2-\theta \sigma^2}.
\end{equation}
\par
\begin{figure}[!ht]
	\begin{center}
	\psfrag{w1}{$\omega_1$}
	\psfrag{w2}{$\omega_2$}
	\psfrag{w3}{$\omega_3$}
	\psfrag{w4}{$\omega_4$}
	\psfrag{$0$}{$0$}
	\psfrag{tan}{ \small{$\tan(\omega T)$}}
	\psfrag{formula}{\small{$\frac{\omega \sigma^2}{\omega^2 \sigma_0^2 + \theta^2 \sigma_0^2 - \theta \sigma^2}$}}
	\psfrag{pi/(2T)}{$\frac{\pi}{2T}$}
	\psfrag{pi/(2T)+pi/T}{$\frac{\pi}{2T} + \frac{\pi}{T}$}
	\psfrag{pi/(2T)+2pi/T}{$\frac{\pi}{2T} + \frac{2\pi}{T}$}
	\psfrag{pi/(2T)+3pi/T}{$\frac{\pi}{2T} + \frac{3\pi}{T}$}
	\includegraphics[width=8cm]{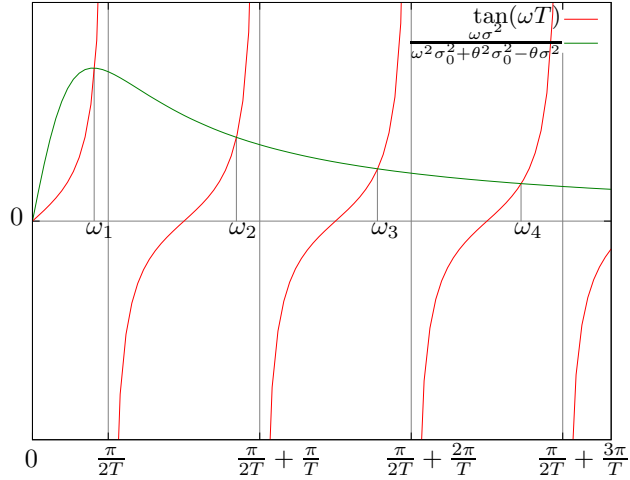}
	\caption{\textbf{(Non-deterministic starting point, $\theta^2 \sigma_0^2-\theta \sigma^2 \geq 0 $)}. Solutions of \eqref{eq:conditional_ornstein_uhlenbeck_functional_equation}. (Ornstein-Uhlenbeck process starting from $X_0 \simdist \mathcal{N}(0,\sigma_0^2)$, $\sigma_0 \neq 0$.) Parameter values are $T=3$, $\sigma=1$, $\theta=3$ and $\sigma_0^2 = 0.4$. }
	\label{fig:inter_ornstein_uhlenbeck_eigen}
	\end{center}
\end{figure}

 \noindent Solutions of this equation are illustrated in Figure \ref{fig:inter_ornstein_uhlenbeck_eigen}. We can show that for any $n \in \N^*$, there is a unique solution of \eqref{eq:general_ornstein_uhlenbeck_functional_equation} in $\left(\frac{n\pi}{T},\frac{n\pi}{T}+\frac{\pi}{2T} \right)$. Moreover a solution lies in $\left(0,\frac{\pi}{2T}\right)$ if and only if $(\theta^2 \sigma_0^2 - \theta \sigma^2)T - \sigma^2 < 0$. 
\par
\begin{figure}[!ht]
	\begin{center}
	\psfrag{w1}{$\omega_1$}
	\psfrag{w2}{$\omega_2$}
	\psfrag{w3}{$\omega_3$}
	\psfrag{w4}{$\omega_4$}
	\psfrag{$0$}{$0$}
	\psfrag{tan}{ \small{$\tan(\omega T)$}}
	\psfrag{formula}{\small{$\frac{\omega \sigma^2}{\omega^2 \sigma_0^2 + \theta^2 \sigma_0^2 - \theta \sigma^2}$}}
	\psfrag{pi/(2T)}{$\frac{\pi}{2T}$}
	\psfrag{pi/(2T)+pi/T}{$\frac{\pi}{2T} + \frac{\pi}{T}$}
	\psfrag{pi/(2T)+2pi/T}{$\frac{\pi}{2T} + \frac{2\pi}{T}$}
	\psfrag{pi/(2T)+3pi/T}{$\frac{\pi}{2T} + \frac{3\pi}{T}$}
	\includegraphics[width=8cm]{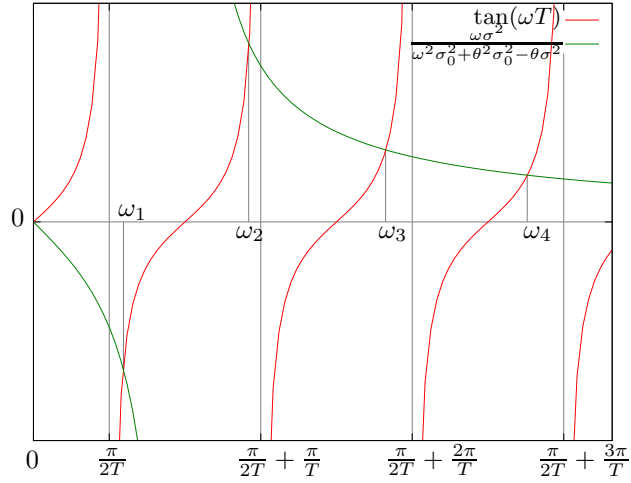}
	\caption{\textbf{(Non-deterministic starting point, $\theta^2 \sigma_0^2-\theta \sigma^2 < 0 $).} Solutions of \eqref{eq:conditional_ornstein_uhlenbeck_functional_equation}. (Ornstein-Uhlenbeck process starting from $X_0 \simdist \mathcal{N}(0,\sigma_0^2)$, $\sigma_0 \neq 0$.) Parameter values are $T=3$, $\sigma=1$, $\theta=3$ and $\sigma_0^2 = 0.3$. }
	\label{fig:general_ornstein_uhlenbeck_eigen}
	\end{center}
\end{figure}

 \noindent \textbf{Second case: $\theta^2 \sigma_0^2-\theta \sigma^2 < 0$.} Here, the term $-\omega^2 \sigma_0^2 + \theta \sigma^2 - \theta^2 \sigma_0^2$ vanishes for $\omega = V:= \sqrt{\theta \frac{\sigma^2}{\sigma_0^2}-\theta^2}$. If $V$ is not a solution of \eqref{eq:ornstein_uhlenbeck_functional_equation}, (\textit{i.e.} if $V$ does not belong to $\left\{ \frac{\pi}{2T} + k \frac{\pi}{T} | k \in \N \right\}$), no other element of this set is a solution, and everything comes again to the same Equation \eqref{eq:general_ornstein_uhlenbeck_functional_equation}. Solutions of this equation are illustrated in Figure \ref{fig:general_ornstein_uhlenbeck_eigen}. We can show that there is a unique solution to \eqref{eq:ornstein_uhlenbeck_functional_equation} in each non-empty interval $\left( \frac{n\pi}{T},\frac{n\pi}{T}+\frac{\pi}{2T} \right) \cap (V, \infty)$ and $\left( \frac{k\pi}{T}-\frac{\pi}{2T}, \frac{k\pi}{T} + \frac{\pi}{2T}\right) \cap \left(0, V\right) $, $k \in \N^*$.
\par In Algorithm \ref{algo:ou_eigen}, we detail the procedure for the computation of the $n$th eigenvalue of the Ornstein-Uhlenbeck covariance operator. The function $\textbf{search}(a,\textnormal{left},\textnormal{right})$ stands for a root finding method. It returns the root of Equation \eqref{eq:ornstein_uhlenbeck_functional_equation} that is bracketed by $[\textnormal{left},\textnormal{right}]$. 
\begin{algorithm}[ht]
	\caption{Ornstein-Uhlenbeck eigenvalue $(\theta,\sigma,\sigma_0,T,n)$}
	\label{algo:ou_eigen}
	\begin{algorithmic}
		\If{$\sigma_0 = 0$} \Comment {\textcolor{Grey1}{There is a unique solution $\omega_n$ of \eqref{eq:ornstein_uhlenbeck_functional_equation} in $\left(\frac{n \pi}{T}-\frac{\pi}{2T},\frac{n \pi}{T}\right)$.}}
			\State \textbf{search}$\left(\omega_n,\frac{n \pi}{T}-\frac{\pi}{2T},\frac{n \pi}{T}\right)$. 
		\vspace{1mm}
		\Else 
		\vspace{1mm}
			\If{$(\theta^2 \sigma_0^2 - \theta \sigma^2) \geq 0$ } \Comment {\textcolor{Grey1}{The vertical asymptote in \eqref{eq:general_ornstein_uhlenbeck_functional_equation} lies on the left of the origin. }}
				\If{$(\theta^2 \sigma_0^2 - \theta \sigma^2)T - \sigma^2 < 0$} \Comment {\textcolor{Grey1}{\eqref{eq:ornstein_uhlenbeck_functional_equation} has a unique solution in $\left(0,\frac{ \pi}{2 T}\right)$.}}
					\State \textbf{search}$\left(\omega_n,\frac{(n-1) \pi}{T},\frac{(n-1) \pi}{T}+\frac{\pi}{2T}\right)$. 
				\Else \Comment {\textcolor{Grey1}{The smallest positive solution $\omega_1$ of \eqref{eq:ornstein_uhlenbeck_functional_equation} lies in $\left(\frac{\pi}{2T},\frac{\pi}{T}\right)$.}}
					\State \textbf{search}$\left(\omega_n,\frac{n \pi}{T},\frac{n \pi}{T}+\frac{\pi}{2T}\right)$. 
				\EndIf
			\Else \Comment {\textcolor{Grey1}{The vertical asymptote of the right-hand side of \eqref{eq:general_ornstein_uhlenbeck_functional_equation} lies on the right the origin. }}
				\If {$\frac{(n-1) \pi}{T} - \frac{\pi}{2T}> \sqrt{\theta \frac{\sigma^2}{\sigma_0^2}-\theta^2}$}
					\State \textbf{search}$\left(\omega_n,\frac{(n-1) \pi}{T},\frac{(n-1) \pi}{T}+\frac{\pi}{2T}\right)$. 
				\ElsIf {$\frac{(n+1) \pi}{T} - \frac{\pi}{2T}< \sqrt{\theta \frac{\sigma^2}{\sigma_0^2}-\theta^2}$}
					\State \textbf{search}$\left(\omega_n,\frac{n \pi}{T} - \frac{\pi}{2T},\frac{n \pi}{T}\right)$. 
				\ElsIf{$\frac{n \pi}{T} - \frac{\pi}{2T}< \sqrt{\theta \frac{\sigma^2}{\sigma_0^2}-\theta^2}$ and $\frac{(n+1) \pi}{T} - \frac{\pi}{2T}> \sqrt{\theta \frac{\sigma^2}{\sigma_0^2}-\theta^2}$}
					\State \textbf{search}$\left(\omega_n,\frac{n \pi}{T} - \frac{\pi}{2T},\sqrt{\theta \frac{\sigma^2}{\sigma_0^2}-\theta^2}\right)$. 
				\Else
					\State \textbf{search}$\left(\omega_n,\sqrt{\theta \frac{\sigma^2}{\sigma_0^2}-\theta^2},\frac{n \pi}{T}-\frac{\pi}{2 T}\right)$. 
				\EndIf
			\EndIf
		\EndIf
		\State \textbf{return} $\lambda_n \gets \frac{\sigma^2}{\omega_n^2+\theta^2}$.
	\end{algorithmic}
\end{algorithm}
\subsubsection{A numerical guess for $\omega_n$.}\label{seq:ornstein_guess}
\par We use $\psi(x) := \frac{\frac{4(8-\pi^2)x^3}{\pi^4}+x}{1-\frac{4 x^2}{\pi^2}}$ as an approximation of $\tan(x)$ on $\left(-\frac{\pi}{2},\frac{\pi}{2}\right)$. We have $\| \tan - \psi \|_{\infty}^{\left(-\frac{\pi}{2},\frac{\pi}{2}\right)} = \frac{10-\pi^2}{2 \pi} \approx 0.02075$. Plugging this into \eqref{eq:conditional_ornstein_uhlenbeck_functional_equation}, we obtain
\begin{equation}\label{eq:conditional_ornstein_uhlenbeck_functional_equation_approx}
\theta \psi(\omega_n T + n\pi) = -\omega_n \hspace{4mm} n \geq 1. 
\end{equation}
This results into a polynomial equation of degree $3$ having a unique (closed-form) solution $\omega^{guess}_n \in \left(\frac{n \pi}{T}-\frac{\pi}{2T},\frac{n \pi}{T}\right)$ which can be used as a starting point for the root finding procedure. 
\section[Closed-form expression for $R_{Y|V}$]{Closed-form expression for $R_{Y|V}$ in the cases of Brownian motion, Brownian bridge and Ornstein-Uhlenbeck processes}\label{sec:annex_R2_computation}
\par We use the same notation as in Section \ref{sec:faster_simulation}. In this Section, we derive closed-form expressions of the matrix $R_{Y|V}:=(\alpha_{ij})_{1\leq i \leq d, 0 \leq j \leq n } \in M_{d,n}(\R)$ corresponding to the affine function $Af_{Y|V}$ defined by $\E[Y|V] = Af_{Y|V}(V)$, in the cases of Brownian motion, Brownian bridge and Ornstein-Uhlenbeck processes.
\par \noindent In the general case, this linear least-square minimization can be performed numerically, but this preliminary stage can become time-consuming when the number of simulation dates grows. 
\par \noindent If $t_0=0 \leq t_1 \leq \cdots \leq t_n=T$ is a subdivision of $[0,T]$, and $X$ is a Gaussian Markov process, we define the affine functions $f^i_j$ by
\begin{equation}\label{eq:closed_form_intermediate_markov}
\E\left[ \int_0^T X_s e_i^X(s) ds \middle| X_{t_0}, \ldots ,X_{t_n} \right] = \sum\limits_{j=0}^{n-1} \E\left[ \int_{t_j}^{t_{j+1}} X_s e_i^X(s) ds \middle| X_{t_j}, X_{t_{j+1}} \right] =: f_j^i(X_{t_j}, X_{t_{j+1}}).
\end{equation}
\subsection{The case of Brownian motion}
\par \noindent Now, assuming that $X = W$ is a standard Brownian motion on $[0,T]$, using Equation \eqref{eq:closed_form_intermediate_markov} we obtain, for $t_j \neq t_{j+1}$, $f_j^i(x,y) = \E\left[ \int_{t_j}^{t_{j+1}} \left( x + \frac{s-t_j}{t_{j+1}-t_j}(y-x) + Y_{s-t_j}^{B,t_{j+1}-t_j} \right) e_i^W(s) ds \right]$, where $Y_{s-t_j}^{B,t_{j+1}-t_j}$ is a standard Brownian bridge on $[t_j,t_{j+1}]$. Hence,
$$
f_j^i(x,y) = x \underbrace{\left( \int_{t_j}^{t_{j+1}} \frac{t_{j+1} - s}{t_{j+1} - t_j} e_i^W(s) ds\right)}_{:=A^i_j} + y \underbrace{\left( \int_{t_j}^{t_{j+1}} \frac{s - t_j}{t_{j+1} - t_j} e_i^W(s) ds\right)}_{:=B^i_j} = x A_j^i + y B_j^i.
$$
\par \noindent Simple algebra leads to
$$
\int_{t_j}^{t_{j+1}} e_i^W(s) ds = \sqrt{\frac{2}{T}} \frac{T}{\pi\left(i-\frac{1}{2}\right)} \left( \cos\left( \pi \left(i-\frac{1}{2}\right) \frac{t_j}{T} \right) - \cos\left(\pi \left(i-\frac{1}{2}\right) \frac{t_{j+1}}{T} \right) \right),
$$
and 
\begin{multline*}
\int_{t_j}^{t_{j+1}} s e_i^W(s) ds = \sqrt{\frac{2}{T}} \frac{T}{\pi\left(i-\frac{1}{2}\right)} \left(t_j\cos\left( \pi \left(i-\frac{1}{2}\right) \frac{t_j}{T} \right) - t_{j+1} \cos\left( \pi \left(i-\frac{1}{2}\right) \frac{t_{j+1}}{T}\right) \right) \\
	+ \sqrt{\frac{2}{T}} \left(\frac{T}{\pi \left(i-\frac{1}{2}\right)}\right)^2 \left(\sin\left( \pi \left(i-\frac{1}{2}\right) \frac{t_{j+1}}{T} \right)-\sin\left( \pi \left(i-\frac{1}{2}\right) \frac{t_j}{T}\right)\right).
\end{multline*}
\par \noindent Hence $\E\left[ \int_0^T W_s e_i^W(s) ds \middle| W_{t_1}, \ldots ,W_{t_n} \right] = \sum\limits_{j=0}^{n-1} A_j^i W_{t_j} + B_j^i W_{t_{j+1}} = \sum\limits_{i=0}^n \alpha_{ij} W_{t_i}$ with, for every $1 \leq j < n$, $\alpha_{ij} = A_j^i + B_{j-1}^i$, $\alpha_{i0} = A^i_0$ and $\alpha_{in} = B^i_{n-1}$. Finally, we get the following closed-form expression for $R_{Y|V}:=(\alpha_{ij})_{1\leq i \leq d, 0 \leq j \leq n }$.
\begin{itemize}
\item If $t_{j-1} < t_j < t_{j+1}$, 
$$
\alpha_{ij} 	= \lambda_i^W \frac{(t_{j+1} - t_{j-1}) e_i^W(t_j) - (t_{j+1} - t_j) e_i^W(t_{j-1})-(t_j-t_{j-1}) e_i^W(t_{j+1})}{(t_{j+1}-t_j)(t_j-t_{j-1})}.
$$
\vspace{3mm}
{\footnotesize
\par If $t_{j-1} = t_j < t_{j+1}$, \ $\alpha_{ij} 	= \lambda_i^W \left(\left(e_i^W\right)'(t_j) - \frac{e_i^W(t_{j+1})- e_i^W(t_j)}{t_{j+1}-t_j} \right).$
\par If $t_{j-1} < t_j = t_{j+1}$, \ $\alpha_{ij} 	= \lambda_i^W \left( \frac{e_i^W(t_j)-e_i^W(t_{j-1})}{t_j-t_{j-1}} -\left(e_i^W\right)'(t_j) \right).$
\par If $t_{j-1} = t_j = t_{j+1}$, \ $\alpha_{ij} = 0.$
}
\item $\alpha_{i0} = 
\left\{\begin{array}{lll}	
\lambda_i^W \left( \left(e_i^W\right)'(t_0) - \frac{e_i^W(t_1)- e_i^W(t_0)}{t_1-t_0} \right) \textnormal{ if } t_1 \neq t_0, \\ 
0 \textnormal{ otherwise.}
\end{array}\right.
$
\item
$\alpha_{in} = 
\left\{\begin{array}{lll}	
\lambda_i^W \left(\frac{e_i^W(t_n)- e_i^W(t_{n-1})}{t_n-t_{n-1}} - \left(e_i^W\right)'(t_n) \right) \textnormal{ if } t_n \neq t_{n-1}, \\ 
0 \textnormal{ otherwise.}
\end{array}\right.
$
\end{itemize}
\par \noindent The equality case is useful when dealing with small time steps that make the numerical evaluation of the divided differences $(e_i^W(t_{j+1})- e_i^W(t_j))/(t_{j+1} - t_j)$ inaccurate.
\subsection{The case of Brownian bridge}
\par \noindent If $X=B$ is a standard Brownian bridge on $[0,T]$, using Equation \eqref{eq:closed_form_intermediate_markov}, we get for $t_j \neq t_{j+1}$, $f_j^i(x,y) = \E\left[ \int_{t_j}^{t_{j+1}} \left( x + \frac{s-t_j}{t_{j+1}-t_j}(y-x) + \left(Y_{s-t_j}^{B,t_{j+1}-t_j}\right) \right) e_i^B(s) ds \right]$, where $Y_{s-t_j}^{B,t_{j+1}-t_j}$ is a standard Brownian bridge on $[t_j,t_{j+1}]$. Hence, very similarly to the case of Brownian motion, 
$$
f_j^i(x,y) = x \underbrace{\left( \int_{t_j}^{t_{j+1}} \frac{t_{j+1} - s}{t_{j+1} - t_j} e_i^B(s) ds\right)}_{:=A^i_j} + y \underbrace{\left( \int_{t_j}^{t_{j+1}} \frac{s - t_j}{t_{j+1} - t_j} e_i^B(s) ds\right)}_{:=B^i_j} = x A_j^i + y B_j^i.
$$
\par \noindent Using that
$$
\int_{t_j}^{t_{j+1}} e_i^B(s) ds = \sqrt{\frac{2}{T}} \frac{T}{\pi i} \left( \cos\left( \pi i \frac{t_j}{T} \right) - \cos\left(\pi i \frac{t_{j+1}}{T} \right) \right),
$$
and
{\small
$$
\int_{t_j}^{t_{j+1}} s e_i^B(s) ds = \sqrt{\frac{2}{T}} \frac{T}{\pi i} \left(t_j\cos\left( \pi i \frac{t_j}{T} \right) - t_{j+1} \cos\left( \pi i \frac{t_{j+1}}{T}\right) \right) + \sqrt{\frac{2}{T}} \left(\frac{T}{\pi i}\right)^2 \left(\sin\left( \pi i \frac{t_{j+1}}{T} \right)-\sin\left( \pi i \frac{t_j}{T}\right) \right),
$$
}
\par \noindent we get $\E\left[ \int_0^T B_s e_i^B(s) ds \middle| B_{t_1},\ldots ,B_{t_n} \right] = \sum\limits_{j=0}^{n-1} \left( A_j^i B_{t_j} + B_j^i B_{t_{j+1}}\right) = \sum\limits_{i=0}^n \alpha_{ij} B_{t_i}$ where, for every $1 \leq j < n$, $\alpha_{ij} = A_j^i + B_{j-1}^i$, $\alpha_{i0} = A^i_0$ and $\alpha_{in} = B^i_{n-1}$. Moreover, we have
$$
A_j^i = \lambda_i^B \left( \left(e_i^B\right)'(t_j) - \frac{e_i^B(t_{j+1}) - e_i^B(t_j)}{t_{j+1} - t_j} \right), \quad \textnormal{and} \quad B_j^i = \lambda_i^B \left( \frac{e_i^B(t_{j+1}) - e_i^B(t_j)}{t_{j+1} - t_j} - \left(e_i^B\right)'(t_{j+1}) \right).
$$
\par \noindent Finally we obtain the following closed-form expression for $R_{Y|V}:=(\alpha_{ij})_{1\leq i \leq d, 0 \leq j \leq n }$.
\begin{itemize}
\item If $t_{j-1} < t_j < t_{j+1}$, 
$$
\alpha_{ij} = \lambda_i^B \frac{(t_{j+1} - t_{j-1}) e_i^B(t_j) - (t_{j+1} - t_j) e_i^B(t_{j-1})-(t_j-t_{j-1}) e_i^B(t_{j+1})}{(t_{j+1}-t_j)(t_j-t_{j-1})}.
$$
\vspace{3mm}
{\footnotesize
\par If $t_{j-1} = t_j < t_{j+1}$, \ $\alpha_{ij} 	= \lambda_i^B \left(\left(e_i^B\right)'(t_j) - \frac{e_i^B(t_{j+1})- e_i^B(t_j)}{t_{j+1}-t_j} \right).$
\par If $t_{j-1} < t_j = t_{j+1}$, \ $\alpha_{ij} 	= \lambda_i^B \left( \frac{e_i^B(t_j)-e_i^B(t_{j-1})}{t_j-t_{j-1}} -\left(e_i^B\right)'(t_j) \right).$
\par If $t_{j-1} = t_j = t_{j+1}$, \ $\alpha_{ij} = 0.$
}
\item $\alpha_{i0} = 
\left\{\begin{array}{lll}	
\lambda_i^B \left( \left(e_i^B\right)'(t_0) - \frac{e_i^B(t_1)- e_i^B(t_0)}{t_1-t_0} \right) \textnormal{ if } t_1 \neq t_0, \\ 
0 \textnormal{ otherwise.}
\end{array}\right.
$
\item $\alpha_{in} = 
\left\{\begin{array}{lll}	
\lambda_i^B \left(\frac{e_i^B(t_n)- e_i^B(t_{n-1})}{t_n-t_{n-1}} - \left(e_i^B\right)'(t_n) \right) \textnormal{ if } t_n \neq t_{n-1}, \\ 
0 \textnormal{ otherwise.}
\end{array}\right.
$
\end{itemize}
\begin{remark}
\par \noindent We obtain the same expression as for Brownian motion, where $(e_n^W,\lambda_n^W)$ is replaced with $(e_n^B,\lambda_n^B)$.
\end{remark}
\subsection{The case of centered Ornstein-Uhlenbeck processes}
\par \noindent If $X$ is an Ornstein-Uhlenbeck process, solution of the SDE $dX_t = -\theta X_t dt + \sigma dW_t$, with $X_0 \simdist \mathcal{N}(0,\sigma_0^2)$ independent of $W$. Consider $t_0=0 \leq t_1 \leq \cdots \leq t_n = T$ a subdivision of $[0,T]$. Using Equation \eqref{eq:closed_form_intermediate_markov} and the conditional Fubini theorem, we obtain
$$
f_i^j(X_{t_j},X_{t_{j+1}}) = \E \left[\int_{t_j}^{t_{j+1}} X_s e_i^{OU}(s)ds \middle| X_{t_j},X_{t_{j+1}} \right] = \int_{t_j}^{t_{j+1}} \E \left[X_s\middle|X_{t_j},X_{t_{j+1}} \right] e_i^{OU}(s)ds,
$$
\par \noindent Assuming that $t_0< t_1 < \cdots t_n $, we easily prove that
$$
\E \left[X_s\middle|X_{t_j},X_{t_{j+1}} \right] = X_{t_j} \frac{e^{\theta \left(t_{j+1}-s\right)}-e^{-\theta \left(t_{j+1}-s\right)}}{e^{\theta \left(t_{j+1}-t_j\right)}-e^{-\theta\left(t_{j+1}-t_j \right)}} + X_{t_{j+1}}\frac{e^{\theta \left(s-t_j\right)} - e^{-\theta \left(s - t_j \right)}}{e^{\theta \left(t_{j+1}-t_j \right)}-e^{-\theta \left(t_{j+1}-t_j \right)}}.
$$
\par \noindent Hence
{
\small
$$
f_i^j(x,y) = x \underbrace{\left( \int_{t_j}^{t_{j+1}} \frac{e^{\theta \left(t_{j+1}-s\right)}-e^{-\theta \left(t_{j+1}-s\right)}}{e^{\theta \left(t_{j+1}-t_j\right)}-e^{-\theta\left(t_{j+1}-t_j \right)}} e_i^{OU}(s)ds \right)}_{:=A_j^i} + y \underbrace{\left( \int_{t_j}^{t_{j+1}} \frac{e^{\theta \left(s-t_j\right)} - e^{-\theta \left(s - t_j \right)}}{e^{\theta \left(t_{j+1}-t_j \right)}-e^{-\theta \left(t_{j+1}-t_j \right)}} e_i^{OU}(s)ds \right)}_{:=B_j^i},
$$
}
where $\left(e_n^{OU}\right)_{n \geq 1}$ are the Karhunen-Loève eigenfunctions of $X$. Finally, we have $\E\left[ \int_0^T X_s e_i^{OU}(s) ds \middle| X_{t_1},\ldots, X_{t_n} \right] = \sum\limits_{j=0}^{n-1} \left( A_j^i X_{t_j} + B_j^i X_{t_{j+1}}\right) = \sum\limits_{i=0}^n \alpha_{ij} X_{t_i}$, with for every $1 \leq j < n$, $\alpha_{ij} = A_j^i + B_{j-1}^i$, $\alpha_{i0} = A^i_0$ and $\alpha_{in} = B^i_{n-1}$. 
\subsubsection{The Ornstein-Uhlenbeck process starting from $0$}\label{eq:starting_from_0}
\par \noindent In this case ($\sigma_0^2 = 0$), as proved in Appendix \ref{sec:karhunen_loeve_basis}, the Karhunen-Loève eigensystem is
\begin{equation}
e_n^{OU}(t) := \frac{1}{\sqrt{\frac{T}{2}-\frac{\sin(2\omega_n T)}{4\omega_n}}} \sin(\omega_n t), \hspace{8mm} \lambda_n^{OU}:= \frac{\sigma^2}{\omega_n^2 + \theta^2}, \hspace{4mm} n \geq 1,
\end{equation}
where $\omega_n$ are the increasingly sorted positive solutions of $\theta \sin(\omega_n T) + \omega_n \cos(\omega_n T) = 0$. This gives 
$$
A_j^i = \frac{-\theta}{\theta^2 + \omega_i^2} \left( \frac{2 e_i^{OU}(t_{j+1})}{e^{\theta (t_{j+1} - t_j)} - e^{-\theta (t_{j+1} - t_j)}} - e_i^{OU}(t_j)\frac{e^{\theta (t_{j+1} - t_j)} + e^{-\theta (t_{j+1} - t_j)}}{e^{\theta (t_{j+1} - t_j)} - e^{-\theta (t_{j+1} - t_j)}} \right) + \frac{1}{\theta^2 + \omega_i^2} \left(e_i^{OU}\right)'(t_j),
$$
\noindent and
$$
B_j^i = \frac{\theta}{\theta^2 + \omega_i^2} \left(e_i^{OU}(t_{j+1})\frac{e^{\theta (t_{j+1} - t_j)} + e^{-\theta (t_{j+1} - t_j)}}{e^{\theta (t_{j+1} - t_j)} - e^{-\theta (t_{j+1} - t_j)}} - \frac{2 e_i^{OU}(t_j)}{e^{\theta (t_{j+1} - t_j)} - e^{-\theta (t_{j+1} - t_j)}} \right) - \frac{1}{\theta^2 + \omega_i^2} \left(e_i^{OU}\right)'(t_{j+1}).
$$
\par \noindent The coefficients $(\alpha_{ij})_{1 \leq i \leq d, 0\leq j \leq n}$ are given by $\alpha_{ij} = A_j^i + B_{j-1}^i$ for $1 \leq j < n$, $\alpha_{i0} = A^i_0$ and $\alpha_{in} = B^i_{n-1}$.
The terms involving $\left(e_i^{OU}\right)'$ vanish. Furthermore, we can show that $\lim\limits_{t_{j+1} \to t_j} A_j^i = 0$ and $\lim\limits_{t_{j-1} \to t_j} B_{j-1}^i = 0$ and deduce the corresponding formula when some dates in the schedule are equal.
\subsubsection{The general Ornstein-Uhlenbeck process}
\par \noindent In this case ($\sigma_0^2>0$), as proved in Appendix \ref{sec:karhunen_loeve_basis}, the Karhunen-Loève eigensystem is given by
\begin{equation}\label{eq:general_KL_OU}
e_n^{OU}(t) := K_n\left(\omega_n \sigma_0^2 \cos(\omega_n t) + \left( \sigma^2 - \theta \sigma_0^2 \right) \sin(\omega_n t)\right), \hspace{8mm} \lambda_n^{OU}:= \frac{\sigma^2}{\omega_n^2 + \theta^2}, \hspace{4mm} n \geq 1,
\end{equation}
where $\omega_n$ are the increasingly sorted positive solutions of $\omega_n \sigma^2 \cos(\omega_n T) + (\theta \sigma^2 - \theta^2 \sigma_0^2 -\omega^2_n \sigma_0^2 ) \sin(\omega_n T) = 0$, and
$$
\frac{1}{K_n^2} = \frac{1}{2 \omega_n} \sigma_0^2 \left(\sigma^2 - \theta \sigma_0^2\right) (1-\cos(2 \omega_n T)) + \frac{1}{2} \sigma_0^4 \omega_n^2 \left(T+\frac{\sin(2 \omega_n T)}{2 \omega_n} \right) + \frac{1}{2} \left(\sigma^2 - \theta \sigma_0^2\right)^2 \left(T-\frac{\sin(2 \omega_n T)}{2 \omega_n} \right).
$$
\par \noindent This gives $e_n^{OU}(t) := K_n \sqrt{\omega_n^2\sigma_0^4 + \left( \sigma^2 - \theta \sigma_0^2 \right)^2} \sin(\omega_n t + \phi_n)$, with $\phi_n = \arccos\left(\frac{\sigma^2 - \theta \sigma_0^2}{\sqrt{\omega_n^2\sigma_0^4 + \left( \sigma^2 - \theta \sigma_0^2 \right)^2}} \right)$ and $\lambda_n^{OU}:= \frac{\sigma^2}{\omega_n^2 + \theta^2}$, $n \geq 1$. Using that for $K \in \R$, $\omega \in \R^*$ and $(t_a,t_b) \in \R^2$, we get
{\small
\begin{multline}
\int_{t_a}^{t_b} \exp(Ks) \sin(\omega s + \phi) ds = \frac{K}{K^2+\omega^2} \left(e^{K t_b} \sin(\omega t_b + \phi) - e^{K t_a} \sin(\omega t_a + \phi) \right) \\
- \frac{\omega}{K^2+\omega^2} \left( e^{K t_b} \cos(\omega t_b + \phi) - e^{K t_a} \cos(\omega t_a + \phi) \right),
\end{multline}
}
we see that the expressions for $(\alpha_{ij})_{1\leq i \leq d, 0 \leq j \leq n}$ established in Section \ref{eq:starting_from_0} remain valid in this case. 

\end{appendices}
\bibliography{biblio}
\end{document}